\documentclass[11pt]{amsart}
\usepackage[pdftex]{graphicx,color}
\usepackage{siunitx}
\def\transp{{\rotatebox[origin=c]{180}{\footnotesize $\perp$}}}
\def\R{\mathbb{R}}
\def\a{\alpha}
\def\b{\beta}
\def\d{\delta}
\def\o{\omega}
\def\k{\kappa}
\let\oldmarginpar\marginpar
\renewcommand\marginpar[1]{
  \oldmarginpar[\raggedleft\footnotesize #1]
  {\raggedright\footnotesize #1}}
\DeclareMathOperator{\diver}{div}
\def\p{\partial}
\newtheorem{algorithm}{Algorithm}

\def\cF{\mathcal{F}}
\def\cT{\mathcal{T}}
\def\V{\mathbb{V}}
\DeclareMathOperator{\trace}{tr}
\def\vrho{\varrho}
\def\g{\gamma}
\def\veps{\varepsilon}
\def\s{\sigma}
\newtheorem{definition}{Definition}
\newtheorem{remark}{Remark}
\numberwithin{remark}{section}

\def\Q{\mathbb{Q}}
\def\cN{\mathcal{N}}
\def\cE{\mathcal{E}}
\def\id{{\rm id}}
\def\tr{{\rm tr}}
\def\W{\mathbb{W}}
\def\G{\mathbb{G}}

\def\meter{{\rm m}}

\def\watt{{\rm W}}

\def\megapascal{{\rm MPa}}
\def\gigapascal{{\rm GPa}}

\graphicspath{{}}
\def\hcI{\widehat{\mathcal{I}}}
\def\bu{\mathbf{u}}
\def\by{\mathbf{y}}
\def\bx{\mathbf{x}}
\def\bq{\mathbf{q}}
\def\bn{\mathbf{n}}
\def\bb{\mathbf{b}}
\def\bw{\mathbf{w}}
\def\bv{\mathbf{v}}
\def\bs{\mathbf{s}}
\def\bz{\mathbf{z}}
\def\be{\mathbf{e}}
\def\bpsi{\pmb{\psi}}
\def\vphi{\varphi}
\def\l{\lambda}
\def\omu{\overline{\mu}}
\def\oa{\overline{\alpha}}
\def\os{\overline{\sigma}}
\def\ok{\overline{\kappa}}
\def\oeta{\overline{\eta}}
\newcommand{\bnu}{{\boldsymbol \nu}}

\def\hby{\widehat{\mathbf{y}}}
\def\hbs{\widehat{\mathbf{s}}}
\def\hbx{\widehat{\mathbf{x}}}
\def\hth{\widehat{\theta}}
\def\hmu{\widehat{\mu}}

\def\htau{\widehat{\tau}}
\def\heps{\widehat{\veps}}
\def\ho{\widehat{\o}}

\newcommand{\modif}[1]{{\color{black}{#1}}}

\newcommand{\rhn}[1]{#1}
\newcommand{\ab}[1]{#1}

\begin{document}
\title[Actuated Bilayers]{Modeling and simulation of thermally
    actuated bilayer plates}
\author[S. Bartels]{S\"oren Bartels}
\address{Department of Applied Mathematics, 
Albert Ludwigs University Freiburg, Germany.}
\email{bartels@mathematik.uni-freiburg.de}
\author[A. Bonito]{Andrea Bonito}
\address{Department of Mathematics, Texas A\&M University, College Station, TX.}
\email{bonito@math.tamu.edu}
\thanks{A.B. was partially supported by NSF Grant DMS-1254618 and AFOSR Grant FA9550-14-1-0234.}
\author[A. H. Muliana]{Anastasia H. Muliana}
\address{Department of Mechanical Engineering, Texas A\&M University, College Station, TX.}
\email{amuliana@tamu.edu}
\thanks{A.H.M. was partially supported by AFOSR Grant FA9550-14-1-0234.}
\author[R. H. Nochetto]{Ricardo H. Nochetto}
\address{Department of Mathematics and Institute for Physical Science
and Technology, University of Maryland, College Park, MD.}
\email{rhn@math.umd.edu}
\thanks{R.H.N was partially supported by NSF Grant DMS-1411808, Institut Henri Poincar\'e and the Simons Visiting Professorship
(Oberwolfach)}
\date{\today}
\subjclass{65N12, 65M60, 35L55, 53C44, 74B20}
\begin{abstract}
We present a mathematical model of polymer bilayers that
undergo large bending deformations when actuated by non-mechanical
stimuli such as thermal effects. The simple model captures a
large class of nonlinear bending effects and can be
discretized with standard plate elements. We devise a 
fully practical iterative scheme and apply it to the simulation
of folding of several practically useful
compliant structures comprising of thin elastic layers.
\end{abstract}

\keywords{Nonlinear elasticity, bilayer bending, finite element method,
iterative solution, micro-switch, deployable structure, encapsulation} 
\maketitle

%%%%%%%%%%%%%%%%%%%%%%%%%%%%%%%%%%%%%%%%%%%%%%%%%%%%%%%%%%%%%%%%%%%%%%%%%%%%%%%%%%
\section{Introduction}
%%%%%%%%%%%%%%%%%%%%%%%%%%%%%%%%%%%%%%%%%%%%%%%%%%%%%%%%%%%%%%%%%%%%%%%%%%%%%%%%%%

\begin{figure}[ht!]
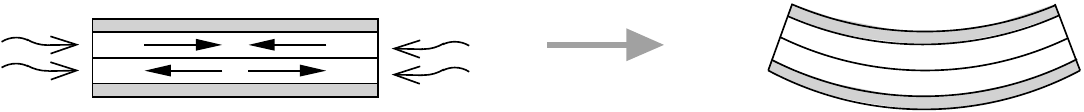
\caption{\label{fig:bilayer_sketch} Schematic description of a thermally actuated
bilayer: heat diffusion starts from the sides into polymer layers that compress and
expand and which are placed between two insulating layers.}
\end{figure}

Bilayer polymers are appealing for the development of autonomous lightweight foldable
structures, such as drug delivery vesicles, flexible constructions for soft robots, 
self-deployable sun
sails in spacecraft, and morphing structures, since they can be manufactured into various shapes
with tunable material properties for each layer. A characteristic feature is
that they can undergo controlled large deformations via relatively small external stimuli. 
Bilayers having two polymeric layers with
different thermo-responsive and swelling characteristics have been shown capable of forming
various practically useful folding shapes \cite{SZTDI12,BALG:10,AOTF93}. 
Different expansion/contraction characteristics of the two layers
when exposed to temperature changes or solvent diffusions lead to out of plane rotations and
curvature changes of the bilayers which trigger folding, cf.~Fig.~\ref{fig:bilayer_sketch}.
When the two layers have significantly different expansion/contraction
behaviors, often characterized by the coefficient of thermal or
moisture expansion, folding can be achieved with relatively smaller stimuli. However, high contrast in mechanical properties, i.e., elastic modulus, between the two layers can induce high stress discontinuities between the layers, leading to delamination. 
Another mechanism for folding and
bending of bilayers can be achieved by integrating two electro-active polymers with the opposite
directions of through-thickness poling axes and when the bilayers are subjected to electric
potential through the thickness one layer would expand while the other would contract 
\cite{TzoTse91,TAJ16,Schm:07a}. One
of the advantages of combining two polymers with different responsive characteristics with
regard to their non-mechanical performances is that the mechanical properties of polymers do
not vary significantly, e.g., extensional elastic moduli of various polymers are typically between
$0.5-5.0\, \gigapascal$, 
which can minimize stress discontinuities at the interfaces between the two layers and
thus can avoid failure due to delamination. Alternatively, bilayers can be formed by combining
two different types of materials, i.e., metal- or ceramic-polymer, which have
significant differences in their mechanical and non-mechanical properties~\cite{SIPL93,KAL08}. 

We are specifically interested in bilayers comprising of two polymeric
layers that can undergo large deformations when exposed to
non-mechanical stimuli, such as temperature changes
or fluid sorption. Due to the slender nature of the bilayers, large
deformations are mainly governed by rotations while the strains in the
bilayer are relatively small~\cite{SRI15,TAJ16}, thereby leading
to negligible stretching and transverse shear effects. 
Furthermore,  the mid-surface of the bilayer plate is considered to be inextensible and non-shearable,
thereby leaving bending as a chief mechanism for shape deformation.
These basic mechanical assumptions lead to a reduced Kirchhoff plate model for the deformation with a \modif{preferred} curvature tensor and an isometry constraint.
The former encodes the mismatch between bilayers while the latter reflects the property that distances among points on the mid surface do not change with shape deformation. In addition, we ignore inertial effects
and assume a quasi-static evolution of the plate.

The bilayer plate is driven by heat conduction or fluid sorption. We
model this diffusion process with a linear Fourier heat conduction or
Fickian diffusion law. Since the deformation is an isometry, the
diffusion equation is insensitive to the plate deformation, whereas
the temperature affects the mismatch between lower and upper layers of
the plate.
In several experimental studies, the bilayers are insulated on top and
bottom and subjected to uniform environmental conditions on their
sides, e.g., at dry conditions the bilayer is immersed in fluid or
from curing at high temperature the bilayer is cooled down to room
temperature or vice versa. In such situations,
the entire polymers will be at uniform temperature or will have uniform fluid content at the steady state and the diffusion process
occurs from the boundaries of the bilayers. 
We model this with homogeneous Neumann boundary condition on top and bottom of the plate,
and either Neumann or Robin boundary conditions on the sides of the
plate. We derive a reduce diffusion model that accounts for these effects.

Our investigations are motivated by lab experiments reporting problems in the controlled
fabrication of nanotubes, cf.~\cite{SZTDI12,Ye130,Yeetal15,Yeetal16}. 
In particular, rectangular bilayer plates occasionally start bending
from the corners leading to formation of so-called {\em dog-ears},
thus failing to attain the desired
cylindrical shape. Our mathematical model and corresponding simulations
indicate that this may be attributed to too rapid changes in the
environment that yield high concentrations
of the diffusing quantity at the corners of the bilayer.
Moreover, we present \modif{computational experiments resembling origami structures}
such as self-assembling cube \cite{SuShMi94,C6MH00195E} and deployable
airfoil \cite{love2007demonstration}, as well as particle encapsulation which is of
interest in drug targeting \cite{StPuIo11}.
Although additional effects become relevant at the nanoscale, see e.g.~\cite{wang2007constitutive}, our computational model could be utilized as a simple tool 
to explore and predict conditions and configurations that enable a controlled 
production of nanoscale and microscale devices.
The method is advantageous with respect to
direct simulation of full three-dimensional thermo-elasticity in terms
of simplicity of implementation, efficiency, performance, and flexibility. 

Our simple two-dimensional model is thus able to capture the principal mechanical
and thermal effects responsible for large deformations observed in lab
experiments of thermally actuated bilayer plates. This suffices for
the aforementioned applications. 
In fact, the present  model and analyses can help designers in simulating desired shape changes and 
determining external stimuli to be prescribed prior to fabricating flexible 
bilayer systems. We refer the reader 
to~\cite{AlBaSm:11,BALG:10,KLPL:05,MaGuSc14,SchEb:01,SIL:95,SIPL93,StPuIo11,SZTDI12,ZLG:14} for further discussion on the applications.
In contrast,
a three-dimensional model capturing strong swelling effects and which
leads to a two-way coupling of the diffusion and deformation processes
has been investigated in~\cite{PAR16}. 

Most of available studies with regards to folding of polymeric
composite structures have been on fabrication and experiments. The
mathematical description and numerical treatment of large bending
deformations of elastic solids has recently undergone some important
development. Dimensionally
reduced models have been derived rigorously from three-dimensional 
hyperelasticity~\cite{FrJaMu02,Schm:07a} and numerical methods capable of approximating
large rotations correctly have been devised and analyzed
in~\cite{Bart13,Bart15-book,BaBoNo15}.
In this work we extend the approach from~\cite{BaBoNo15} by including temperature
dependence in the \modif{preferred} curvature tensor. The mechanical part
consists of a nonlinear
Kirchhoff model for which we devise a finite element discretization based on standard
plate elements. 
The heat equation decouples from the mechanical equation due to the
inextensibility of the plate and we approximate it with standard
finite element methods together with the backward Euler time stepping algorithm.
The proposed iterative numerical method for the coupled system
is \modif{roughly ten times faster} than the one used in~\cite{BaBoNo15}. 

The outline of this article is as follows. In Section~\ref{sec:math_model}
we introduce the thermo-mechanical mathematical model and describe its
dimension reduction. Corresponding partial differential equations
are formulated in Section~\ref{sec:pdes}. We then present in
Section~\ref{sec:num_scheme} the temporal and spatial discretization of
the nonlinear and constrained system composed of second order diffusion
and fourth order bending equations. We report several intruiguing
and practically useful numerical experiments
in Section~\ref{sec:experiments}.

%%%%%%%%%%%%%%%%%%%%%%%%%%%%%%%%%%%%%%%%%%%%%%%%%%%%%%%%%%%%%%%%%%%%%%%%%%%%%%%%%%
\section{Mathematical model}\label{sec:math_model}
%%%%%%%%%%%%%%%%%%%%%%%%%%%%%%%%%%%%%%%%%%%%%%%%%%%%%%%%%%%%%%%%%%%%%%%%%%%%%%%%%%

%---------------------------------------------------------------------------------
\subsection{Hyperelastic materials}
%---------------------------------------------------------------------------------
We model polymers as isotropic and elastic materials.
According to St. Venant--Kirchhoff description, 
the mechanical behavior of the system is governed by a hyperelastic
stored energy density
\[
W(E) := \mu |E|^2 + \frac{\l}{2}  (\tr E)^2,
\]
where $E := F^\transp F-I$ is the Green--Lagrange strain tensor related to the 
deformation gradient $F$, 
and $\lambda,\mu$ are the (first and second) Lam\'e constants
\cite{Ogde84-book}; hereafter we let $|E|^2 = E:E$ be the square of the 
Frobenius norm of $E$ and $\tr E = E:I$ be the trace of $E$.
As strain tensor we use the temperature dependent
quantity $E_\Theta := F^\transp F-(1+\a \Theta) I \in \R^{3\times 3}$ which generalizes the Green--Lagrange 
strain tensor. Here $\a, \Theta $ are the thermal expansion
coefficient and temperature change. Therefore, a simple calculation
yields
\[
W(E_\Theta) = W(E) - (2\mu+3\lambda) \alpha\Theta \tr E + c_\Theta
\]
with $ c_\Theta := (3\mu+\frac{9}{2}\lambda)\alpha^2 \Theta^2$.
We note that the elastic material constants $\lambda,\mu$
can also change with temperature or fluid
concentration. 
We simplify the model further upon realizing that
$(\tr E)^2 \le 3|E|^2$, whence
\[
\mu |E_\Theta|^2 \le W(E_\Theta) \le \big( \mu + \frac{3 \l}{2}\big) |E_\Theta|^2.
\]
This equivalence motivates our choice of energy
\begin{equation}\label{simplified-energy}
  W_\Theta(F) := \frac{\mu}{4} |E_\Theta|^2
  = \frac{\mu}{4} \big|F^\transp F-(1+\a \Theta) I \big|^2.
\end{equation}
which in turn allows us to follow the arguments of~\cite{BaBoNo15} for the
formal dimension reduction of the plate model. A rigorous derivation for more general
material models including the above St. Venant--Kirchhoff  material can be
found in~\cite{Schm:07a} which leads to the same dimensionally reduced
model up to a different prefactor in the energy density. 
For further details on admissible energies and relaxation effects of the dimension reduction,
we refer the reader to \cite{FrJaMu:06,FrJaMu02,PAR16,Schm:07a}. We remark that a more general,
nonlinear dependence of the strain tensor on temperature can be used;
for our experiments a linear relation turned out to be sufficient.
%----------------------------------------------------------------------------------
\subsection{Bilayer configuration and evolution hypotheses}
%----------------------------------------------------------------------------------

We consider the following configuration below. We let
$\omega \subset \mathbb{R}^2$ be the flat parametric domain, $\delta>0$ be the
thickness of the plate, and $\omega_\delta := \omega \times (-\delta/2, \delta/2)$
be the plate in the undeformed configuration. We let $\bx := (\bx',x_3) \in\omega_\delta$
denote a generic point and $t\in (0,\infty)$ time. We further indicate with
$\bu (\cdot,t): \omega_\delta\to\mathbb{R}^3$ the deformation of the
plate and $\gamma_\delta(t) = \bu(\cdot,t)(\omega_\delta)$ 
the deformed configuration of the plate at time $t$.

For the evolution model considered below we
assume that inertial effects are insignificant and that no mechanical 
dissipation occurs so that the plate immediately
adjusts to temperature changes. This yields a dynamics driven by 
temperature, which mathematically entails that $\Theta$ solves a suitable diffusion
equation on the flexible surface and $\bu(\cdot,t)$ is a minimizer of an elastic energy
associated with the density $W_\Theta(\nabla \bu(\cdot,t))$.
Due to the assumptions of thermal insulation
of the bilayer top and bottom and of mid-surface inextensibility,
it turns out that the diffusion equation decouples from the deformation equation. 
We examine these two models below and outline reduced models for the
vanishing thickness limit $\delta\to0$.

%--------------------------------------------------------------------------------
\subsection{Reduced bilayer plate model}
%--------------------------------------------------------------------------------
%
We adjust the simplified energy \eqref{simplified-energy} for hyperelastic
materials to model thin bilayers and derive effective
energies describing large deformations. In particular, we consider two layers
of materials glued
on top of each other with {\em different} thermal material constants
such as the thermal 
expansion coefficients. We assume that one material expands and the other
compresses above a critical temperature, which we assume to be zero
for simplicity. Therefore, we consider \eqref{simplified-energy}
with $F=\nabla\bu$
\[
W_\Theta(\nabla\bu,\bx,t) = \frac{\mu(\bx)}{4}
\Big| \nabla\bu^\transp \nabla\bu - (1+\a(\bx)
\Theta(\bx,t)) I \Big|^2,
\]
where the coefficients $\a(\bx), \mu(\bx)$ jump across $x_3=0$
\[
\a(\bx) := \begin{cases}
+\a & \mbox{for } x_3 >0, \\
-\a & \mbox{for } x_3 <0,
\end{cases} \qquad
\mu(\bx) := \begin{cases}
\mu^+(x') & \mbox{for } x_3 >0, \\
\mu^-(x') & \mbox{for } x_3 <0.
\end{cases}
\]
We follow~\cite{FrJaMu:06,Schm:07b,BaBoNo15}
to identify a dimensionally reduced model corresponding to the
limit $\d \to 0$. Formally, this is based on the Kirchhoff assumption that
the actual deformation $\bu$, subject to given
forces and boundary conditions can be, up to higher order contributions,
represented as 
\begin{equation}\label{eq:ass_kirchhoff}
\bu(\bx',x_3,t) = \by(\bx',t) + x_3 \bb(\bx',t)
\end{equation}
with a mapping $\by: \o\times(0,\infty) \to \R^3$ that describes the
deformation of the midplane
$\o\times \{0\}$ and a vector field $\bb: \o \times (0,\infty) \to \R^3$ that
is normal to the deformed midplane $\gamma(t)=\by(\o,t)$.
This means that fibers perpendicular to $\omega$ in the undeformed
configuration remain normal to the mid-surface $\gamma(t)$ in the
deformed configuration. For ease of presentation we do not write the argument 
$t$ explicitly in the remainder of this subsection.
Inserting the corresponding deformation
gradient
\[
\nabla \bu (\bx) = [\nabla' \by(\bx'),\bb(\bx')] + x_3 [\nabla'\bb(\bx'),0]
\]
into the scaled elastic energy functional
\[
I_\Theta[\bu]= \frac{1}{\d^3} \int_{\o_\d} W_\Theta (\nabla\bu,\cdot)
\]
leads to
\[
\begin{split}
I_\Theta [\bu] = \frac{1}{\d^3} \int_\o \int_{-\d/2}^{\d/2}
\frac{\mu}{4}\Bigg\{\Bigg|
& \begin{bmatrix} (\nabla' \by)^\transp \nabla \by - (1 \pm \a
\Theta) I_2 & 0 \\ 0 & |\bb|^2 - (1\pm \a \Theta)
\end{bmatrix} \\
+ x_3
& \begin{bmatrix} (\nabla' \bb)^\transp \nabla' \by +
(\nabla' \by)^\transp \nabla' \bb & (\nabla' \bb)^\transp \bb \\
\bb^\transp \nabla' \bb & 0 \end{bmatrix} 
\\
+ x_3^2
& \begin{bmatrix} (\nabla' \bb)^\transp \nabla' \bb & 0
\\ 0 & 0 \end{bmatrix}
\Bigg|^2 \Bigg\}.
\end{split}
\]
The scaling of the elastic energy by $\d^{-3}$ corresponds to deformations 
that describe a bending behavior of the thin plate and is
crucial for identifying 
limiting equations. We assume that $I_\Theta[\bu]$ remains
bounded as $\d\to 0$,
carry out the integration in $x_3$ direction, and deduce
necessary scaling properties
of terms arising in the energy functional; details can be found in
\cite{BaBoNo15}.
The first necessary condition to have a finite limit as $\d\to 0$
is that $\bb$ has length $|\bb|=1\pm\alpha\Theta$ in the upper 
and lower layers. Assuming that \modif{$\a/\d$ is finite}, which is explained
below, shows that for $\d \to 0$ we have that $\bb$ equals the unit
normal $\bnu$ to the deformed midplane $\gamma$ and that
$(\nabla'\bb)^\transp \bb = 0$. 
Let $G$ and $H$ be the first and second fundamental forms of $\gamma$,
i.e., the symmetric matrices
\[
G = (\nabla' \by)^\transp \nabla' \by, \quad H = - (\nabla'
\bnu)^\transp \nabla' \by.
\]
Then we find the second necessary condition  $G = I_2$
for a finite limit $\d\to 0$,
\begin{equation*}
[\nabla' \by]^T \nabla' \by = I_2,
\end{equation*}
i.e., $\by: \o \to \R^3$ 
defines an isometric deformation of the midplane $\omega$. With this we derive
the central identity
\[
\d^{-3} \int_{-\d/2}^{\d/2} \frac{\mu}{4}  |2 x_3 H \mp \a \Theta I_2 |^2  
= \frac{\omu}{12}  \Big( |H - \oa  \theta I_2|^2  - \frac43 \oa^2 \theta^2 \Big),
\]
in which $\omu$ is the average of $\mu(\cdot,\pm x_3)$, 
$\theta(\bx') = \Theta(\bx',0)$ is the effective temperature that obeys a reduced
diffusion equation discussed below, and $\oa$ is the 
effective thermal expansion coefficient per unit thickness given by
\begin{equation}\label{d:baralph}
\oa := \lim_{\d\to 0} \frac{3\,\a}{\d}.
\end{equation}
\modif{Note that $\overline{\a} \theta$
has the unit of a curvature, i.e. \SI{}{\per \meter}.}

We thus  infer that the reduced energy functional
takes the form 
\begin{equation}\label{reduced-energy}
E[\by] = \frac{1}{12} \int_\o \omu \big|H- \oa \theta I_2\big|^2,
\end{equation}
up to constant terms that we ignore from now on.
Assuming that $\a$ is comparable with the plate thickness
is necessary to avoid delamination and obtain
a total bending effect of the thin bilayer plate resulting
from a compressive and expansive behavior of the individual layers.
This also shows that small material differences can lead to large deformations.

%--------------------------------------------------------------------------------
\subsection{Reduced diffusion model}
%--------------------------------------------------------------------------------
%
As in the mechanical part we allow the two different materials in the upper and
lower layers of the plate to possess different physical constants, namely
 heat capacities $\sigma^{\pm}$ and
conductivities $\k^{\pm}$. We define
\begin{equation}\label{meanvalues}
\os := \frac{1}{2} (\sigma^+ + \sigma^-),
\quad
\ok := \frac{1}{2} (\k^+ + \k^-)
\end{equation}
to be their mean values in the transversal direction.
We make the key simplifying assumptions that the plate is
thermally insulated on top and bottom and that
diffusion is taking place mostly in tangential direction of
the plate. This
means that the gradient of temperature, and so the heat flux $\bq=
-\k\nabla\Theta$ which obeys the Fourier law of diffusion, satisfies
the orthogonal decomposition
\begin{equation}\label{temperature}
\nabla\Theta = \nabla_\gamma\Theta + \partial_\bnu\Theta \, \bnu,
\end{equation}
where for $(\bx',x_3) \in \omega_\d$
\begin{gather*}
\nabla_\gamma \Theta(\bx',x_3) = \big(I-\bnu(\bx')\otimes\bnu(\bx')\big)
\nabla\Theta(\bx',x_3),
\\
\partial_\bnu \Theta(\bx',x_3) = \bnu(\bx') \cdot \nabla\Theta(\bx',x_3),
\end{gather*}
as well as
\begin{equation}\label{normal-flux}
\partial_\bnu \Theta(\bx',x_3) \to 0 \quad\text{as } \d \to 0.
\end{equation}
We start with the energy balance in the slender 3d set $\gamma_\delta(t) :=
\bu(\cdot,t)(\omega_\delta)$, i.e., in the
deformed plate with positive thickness $\d$ where $\bu$ is defined
in~\eqref{eq:ass_kirchhoff}:
\begin{equation}\label{energy-balance}
\frac{d}{dt} \int_{\gamma_\delta(t)} \sigma \Theta
= - \int_{\partial\gamma_\delta(t)} \bq\cdot\bn
= - \int_{\gamma_\delta(t)} \diver\bq.
\end{equation}
Using Reynolds' transport theorem, cf., e.g.,~\cite{DziukElliott07},
we can write the left-hand side as follows:
\[
\frac{d}{dt} \int_{\gamma_\delta(t)} \sigma \Theta
= \int_{\gamma_\delta(t)} \sigma \big(\partial_t\Theta
+ \nabla\Theta \cdot\partial_t\bu + \Theta \diver \partial_t\bu \big).
\]
Using \eqref{temperature}, we obtain for the second summand
in the right-hand side as $\delta\to0$
\[
\frac{1}{\delta} \int_{\gamma_\delta(t)} \sigma \nabla \Theta \cdot
\partial_t\bu
= \frac{1}{\delta} \int_{\gamma_\delta(t)} \sigma\nabla_\gamma \Theta \cdot
\partial_t\bu + \sigma \partial_\bnu \Theta \, \bnu\cdot\partial_t\bu
\to \int_{\gamma(t)} \os \nabla_\gamma \theta\cdot\partial_t\by,
\]
where $\theta(\bx',t) = \Theta(\bx',0,t)$ is the temperature in the
midsurface $\gamma(t)$ and $\os$ is defined in \eqref{meanvalues}.
It remains to evaluate the term $\diver \partial_t\bu$.

To do so, we set $\mathbf{F}(\bx)=\mathbf{f}(\bu)$
for a vector-valued function $\mathbf{f}:\gamma_\d\to\R$ and use
\eqref{eq:ass_kirchhoff} to deduce
\[
\partial_{x_3} \mathbf{F}(\bx) = \nabla \mathbf{f}(\bu(\bx))
\partial_{x_3} \bu(\bx) =
\nabla \mathbf{f} (\bu(\bx)) \bb(\bx')
\qquad\forall \, \bx\in\gamma_\d.
\]
Since the divergence is the trace of the gradient, we invoke the
decomposition \eqref{temperature} for each component of $\mathbf{f}$
together with the preceding expression to obtain
\begin{equation}\label{divergence}
\diver \mathbf{f} = \diver_\gamma \mathbf{f} + \bnu\cdot\nabla\mathbf{f} \, \bnu
= \diver_\gamma \mathbf{f} + |\bb|^{-2} \bb\cdot\partial_{x_3}\mathbf{F}.
\end{equation}
To apply this formula to $\mathbf{f}=\partial_t\bu
= \partial_t \by + x_3 \, \partial_t \bb$ we realize that
$\partial_{x_3}\partial_t\bu = \partial_t \bb$ whence
$\bb\cdot\partial_{x_3}\partial_t\bu = \bb\cdot\partial_t\bb
= \frac12\partial_t|\bb|^2=0$ and
\[
\diver \partial_t\bu = \diver_\gamma \partial_t\bu.
\]
Therefore, as $\delta\to0$ we arrive at
\[
\frac{d}{dt} \frac{1}{\delta} \int_{\gamma_\delta(t)} \sigma \Theta
\to \int_{\gamma} \os \big(\partial_t\theta
+ \nabla_\gamma \theta \cdot\partial_t\by + \theta \diver_\gamma
\partial_t\by \big)
= \frac{d}{dt} \int_{\gamma(t)} \os \theta.
\]
We now deal with the flux term $\bq = - \k \nabla \Theta$. We first write
\[
\diver \bq = \diver_\gamma \bq + \diver_\gamma^\perp \bq,
\]
where $\diver_\gamma^\perp \bq = \bnu \cdot \nabla\bq \, \bnu$ in
view of \eqref{divergence}.
Integrating by parts in the normal direction $\bnu$, and using the vanishing
Neumman boundary conditions $\bq^+\cdot\bnu=\bq^-\cdot\bnu=0$ on top
and bottom of the plate $\gamma_\d$ gives
\begin{equation}\label{div-perp}
\int_{\gamma_\delta} \diver_\gamma^\perp \bq = 0.
\end{equation}
On the other hand, in light of \eqref{temperature} and \eqref{normal-flux}
we deduce as $\delta\to0$
\begin{equation*}
\frac{1}{\delta} \int_{\gamma_\delta} \diver_\gamma \k \, \nabla\Theta
= \frac{1}{\delta} \int_{\gamma_\delta} \diver_\gamma \k \,\nabla\gamma\Theta
+ \frac{1}{\delta} \int_{\gamma_\delta} \diver_\gamma \k \,
(\partial_\bnu\Theta \, \bnu )
\to \int_{\gamma(t)}\diver_\gamma \ok \nabla_\gamma\theta,
\end{equation*}
where $\ok$ is defined in \eqref{meanvalues} because
\[
\diver_\gamma (\partial_\bnu\Theta \, \bnu )
= \nabla_\gamma(\partial_\bnu\Theta)\cdot\bnu + \partial_\bnu\Theta
\diver_\gamma\bnu = \partial_\bnu\Theta \, \tr H \to 0
\quad\text{as } \d\to0.
\]

Collecting all the previous results, we conclude that the limit of
\eqref{energy-balance} as $\delta\to0$ is the surface conservation equation
\begin{equation}\label{heat-reduced}
  \frac{d}{dt} \int_{\gamma(t)} \os \theta -
  \int_{\gamma(t)} \diver_\gamma (\ok \nabla_\gamma \theta) = 0.
\end{equation}  
This equation is consistent with
the diffusion equation on a surface derived in~\cite{DziukElliott07}.

We conclude with a simple extension of \eqref{heat-reduced} which accounts for diffusion {\it transversal}
to the plate. Suppose that the normal fluxes $\bq^+,\bq^-$ on top and
bottom of the plate $\gamma_\d(t)$ do not vanish but rather scale
proportional to $\d$. In this case, \eqref{div-perp} reduces to
\[
\frac{1}{\delta} \int_{\gamma_\d} \diver_\gamma^\perp\bq =
\frac{\bq^+\cdot\bnu - \bq^-\cdot\bnu}{\delta} \to \int_\gamma f
\qquad\text{as } \d\to0.
\]
In other works, the function $f$ acts as an effective source term in
the surface diffusion equation. In addition, the following limit does
no longer vanish
\[
\frac{1}{\delta} \int_{\gamma-d} \sigma \partial_\bnu\Theta \,\bnu\cdot
\partial_t\bu + \k \partial_\bnu\Theta \, \tr H
\to \int_{\gamma_d} \partial_\bnu\Theta \big(\bar{\sigma} V + \bar\k \, \tr H\big)
\qquad\text{as }\d\to0,
\]
where $V=\bnu\cdot\partial_t\by$ is the normal velocity
and $h=-\tr H$ is the mean curvature  of $\gamma(t)$. This leads to
the following variant of \eqref{heat-reduced}
\[
\frac{d}{dt} \int_{\gamma(t)} \bar\sigma \theta +
  \int_{\gamma(t)} - \diver_\gamma(\bar\k \nabla_\gamma\theta)
  + \partial_\bnu\theta \big(\bar{\sigma} V - \bar\k \, h\big)
  + f = 0,
\]
which of course must be supplemented with a diffusion equation in the
surroundings of $\gamma(t)$ to determine the quantity
$\partial_\bnu\theta$. In contrast to \eqref{heat-reduced}, this form
of the PDE does couple diffusion and plate geometry.

%%%%%%%%%%%%%%%%%%%%%%%%%%%%%%%%%%%%%%%%%%%%%%%%%%%%%%%%%%%%%%%%%%%%%%%%%%%%%%%%%%
\section{Governing PDEs: Weak forms}\label{sec:pdes}
%%%%%%%%%%%%%%%%%%%%%%%%%%%%%%%%%%%%%%%%%%%%%%%%%%%%%%%%%%%%%%%%%%%%%%%%%%%%%%%%%%

The governing equations are formulated in the cylinder $\omega\times(0,T)$ and
the independent variables are denoted by $(\bx,t) \in \omega\times(0,T)$
for simplicity.

%---------------------------------------------------------------------------------
\subsection{Plate equation}
%---------------------------------------------------------------------------------
We first derive the Euler-Lagrange equation for a minimizer of the
bending energy \eqref{reduced-energy} given a fixed temperature distribution.

\ab{In \eqref{reduced-energy}, we omit writing the variable $\bx$ and write
\[
E[\by] = \frac{1}{12} \int_\o \omu \big| H - \oa\theta I_2 \big|^2,
\]
}
For isometries, the $ij$-th element of $H$, namely
$h_{ij} = \partial_i\partial_j \by \cdot \bnu$, satisfies the key relation
\begin{equation}\label{sec-fund}
\partial_i\partial_j \by = h_{ij} \bnu;
\end{equation}
hence $\partial_i\partial_j \by$ is parallel to $\bnu$. This
immediately implies equality of the Frobenius norm of the
second fundamental form $H$ and the Hessian of $\by$
\[
|H| = |D^2\by|.
\]
Developing the square yields the following equivalent expression for
$E[\by]$
\[
E[\by] = \ab{\frac{1}{12}} \int_\o \omu |D^2\by|^2 - 2\oa\omu\theta \, H:I_2
+ 2 \omu (\oa\theta)^2.
\]
Since $\bnu = \partial_1\by\times \partial_2\by$ for isometries,
in view of \eqref{sec-fund} we deduce the expression
\[
H:I_2 = \Delta\by \cdot \big(\partial_1\by\times \partial_2\by\big).
\]
Consequently, $\by$ is a minimizer of the energy
\[
E[\by] = \ab{\frac{1}{12}} \int_\o \omu |D^2\by|^2 - 2\oa\,\omu\theta \,
\Delta\by \cdot \big( \partial_1\by\times \partial_2\by\big)
+ 2 \omu (\oa\theta)^2
\]
subject to the isometry constraint
\begin{equation}\label{isometry}
[\nabla \by]^\transp \nabla \by = I_2.
\end{equation}
Computing the first variation of these two relations with respect to
$\by$ yields the weak form of the Euler-Lagrange equation
\begin{equation}\label{EL-1}
\begin{aligned}
(\omu D^2\by,D^2\bw) &- (\Delta \bw \cdot [\p_1 \by \times \p_2 \by],\oa \omu \theta) 
\\
& - (\Delta \by \cdot[\p_1 \bw \times \p_2 \by + \p_1 \by \times \p_2 \bw], \oa \omu \theta)
=0
\end{aligned}
\end{equation}
for any test function $\bw$ satisfying the linearized isometry condition
\begin{equation}\label{iso}
[\nabla \bw]^\transp \nabla \by + [\nabla \by]^\transp \nabla \bw = 0.
\end{equation}
The identity \eqref{EL-1} can be simplified further because the
third term on the left-hand side vanishes. To see this, we first
recall from \eqref{sec-fund} that second derivatives
$\p_i\p_j\by$ of $\by$ are parallel to $\bnu = \p_1\by\times\p_2\by$
for isometries, whence the mean curvature $h = -\tr (H G^{-1}) = -\tr H$ satisfies
\[
- \Delta \by = h \, \p_1\by\times\p_2\by.
\]
We next observe that $\partial_i\by \cdot \partial_j\by=0$ for
$i\ne j$ according to \eqref{isometry} and
$\partial_i\by \cdot \partial_i\bw = 0$ because of \eqref{iso}.
Using now the formula
$(a\times b)\cdot (c\times d) = (a\cdot c)(b\cdot d) - (a\cdot d)(b\cdot d)$ yields 
\begin{align*}
-\Delta \by & \cdot [\p_1 \bw \times \p_2 \by + \p_1 \by \times \p_2 \bw] \\
&= h (\p_1 \by \cdot \p_1 \bw) (\p_2 \by \cdot \p_2 \by) 
- h (\p_1 \by \cdot \p_2 \by) (\p_2 \by \cdot \p_1 \bw) \\
& \qquad + h (\p_1 \by \cdot \p_1 \by) (\p_2 \by \cdot \p_2 \bw) 
- h (\p_1 \by \cdot \p_2 \bw) (\p_2 \by \cdot \p_1 \by) =0.
\end{align*}
This implies that $\by$ is a solution of the simplified fourth order equation
\begin{equation*}
(\omu D^2\by,D^2\bw) - (\Delta \bw \cdot [\p_1 \by \times \p_2 \by],\oa
  \omu \theta) = 0
\end{equation*}
for all test functions $\bw$ satisfying the linearized isometry
condition \eqref{iso}.

%--------------------------------------------------------------------------------
\subsection{Diffusion equation}
%--------------------------------------------------------------------------------

We now turn to the diffusion equation \eqref{heat-reduced}. Since its derivation 
is valid for any patch in $\gamma(t)$, we realize that it implies the strong form 
\begin{equation}\label{e:diff_eq}
\os D_t \theta + \os \theta \diver_\gamma\p_t\by - \diver_\gamma
(\ok\nabla_\gamma\theta) = 0,
\end{equation}
where $D_t\theta = \p_t \theta + \nabla_\gamma \theta \cdot\p_t\by$ is
the material derivative of $\theta$, cf.~\cite{DziukElliott07}. 
We split the boundary $\partial \gamma(t)$ of $\gamma(t)$ in two disjoint pieces  
$\partial_D \gamma(t)$ and $\partial_R \gamma(t)$, where $\partial_D \gamma(t)$ 
is the portion of $\partial \gamma(t)$ where we prescribe the temperature and  
$\partial_R\gamma(t)$ is where we impose the Robin condition 
\[
\ok \nabla_\gamma \theta \cdot\bn = \oeta (\theta_{ext} - \theta).
\]
 
Multiplying \eqref{e:diff_eq} by any test function $\vphi$ with vanishing material derivative $D_t\vphi=0$ and vanishing on $\partial_D \gamma(t)$, and integrating over
$\gamma(t)$ leads to 
\[
\int_{\gamma(t)} \os D_t (\theta\vphi) + \os \theta\vphi \diver_\gamma\p_t\by + 
\ok\nabla_\gamma\theta \cdot \nabla_\gamma \vphi =
\int_{\partial_R\gamma(t)} \ok \nabla_\g \theta \cdot\bn \, \vphi.
\]
We thus have that 
\begin{equation}\label{heat-weak-1}
  \frac{d}{dt} \int_{\gamma(t)} \os \theta \vphi
  +  \int_{\gamma(t)} \ok \nabla_\gamma\theta \cdot \nabla_\gamma \vphi
   + \int_{\partial_R \gamma(t)} \oeta  \theta \vphi 
  = \int_{\partial_R\gamma(t)} \oeta  \theta_{ext} \vphi;
\end{equation}
this form of the diffusion equation on $\gamma(t)$ is due to~\cite{DziukElliott07}.
In the present context, this equation
simplifies further because $\by$ is an isometry for all $t$. In fact,
since the first fundamental form satisfies $G=[\nabla \by]^\transp \nabla \by= I_2$ we can write with
$\bz = \by(\bx)$ 
\[
\int_{\gamma(t)} \os(\by^{-1}(\bz))\theta(\bz,t)\vphi(\bz,t) d\bz
= \int_\omega \os(\bx)\theta(\by(\bx,t),t)\vphi(\by(\bx,t),t) d\bx,
\]
and
\begin{align*}
\int_{\gamma(t)} \ok (\by^{-1}(\bz))\nabla_\gamma
  \theta(\bz,t) & \cdot\nabla_\gamma\vphi(\bz,t) d\bz
\\
&= \int_\omega \ok(\bx) \nabla \theta(\by(\bx,t),t)\cdot\nabla\vphi(\by(\bx,t),t) d\bx
\end{align*}
because $\nabla_\gamma\vphi(\bz,\cdot)=\nabla\by G^{-1} \nabla
\vphi(\by(\bx,\cdot),\cdot)$.
This allows us to express \eqref{heat-weak-1} in the parametric domain
$\omega$ and avoid the dependence on $t$ in the test function
$\vphi$. We thus get the following simple weak form of the diffusion
equation
\begin{equation*}
\int_\o \os \partial_t \theta(\cdot,t) \, \vphi + \ok \nabla\theta(\cdot,t) \cdot
\nabla\vphi + \int_{\partial_R \omega} \oeta \theta \vphi 
= \int_{\partial_R\omega} \oeta \theta_{ext} \vphi, 
\end{equation*}
i.e., it is sufficient to solve the diffusion equation in the reference
configuration $\omega$. Moreover, the diffusion equation decouples
from the plate equation due to the isometry property and the vanishing Neumann
boundary conditions for temperature assumed on top and bottom of the plate. 

%%%%%%%%%%%%%%%%%%%%%%%%%%%%%%%%%%%%%%%%%%%%%%%%%%%%%%%%%%%%%%%%%%%%%%%%%%%%%%%%%%
\section{Numerical scheme}\label{sec:num_scheme}
%%%%%%%%%%%%%%%%%%%%%%%%%%%%%%%%%%%%%%%%%%%%%%%%%%%%%%%%%%%%%%%%%%%%%%%%%%%%%%%%%%

%---------------------------------------------------------------------------------
\subsection{Time discretization}
%---------------------------------------------------------------------------------

We next describe a discrete-time scheme to compute the evolution
of thermally induced bilayer bending effects. We consider clamped boundary
conditions for the mechanical equation and Dirichlet conditions for the 
diffusion process on the subset $\p_D\o \subset \p \o$. 
A Robin-type boundary condition is imposed on the remaining part
$\p_R \o = \p\o \setminus \p_D \o$ while no explicit boundary conditions
are imposed on the deformation on this part. 
We further impose the condition
that deformations are contained in a convex set $K$ which models
the presence of an obstacle, e.g.,
$K= \{ \by=(y_i)_{i=1}^3 \in L^2(\o;\R^3): y_3 \le 1 \}$,
or $K$ is the entire space $K=L^2(\o;\R^3)$ in the case of no obstacle.
The $L^2$ scalar product on $\o$ of functions or vector fields is
denoted by $(\cdot,\cdot)$. If the $L^2$ scalar product is taken on a 
set $A$ we write $(\cdot,\cdot)_A$.
In addition, we let $\{\be_1,\be_2\}$ be the canonical unit
vectors in $\omega$ and use the notation $d_t\theta^{k+1}$ to indicate
the scaled backward difference
\[
d_t\theta^{k+1} := \tau^{-1} \big( \theta^{k+1} - \theta^k \big).
\]

\begin{algorithm}[abstract time stepping]\label{alg:alg_abstract}
Let $\tau>0$ be a uniform time step, set $\theta^0 :=0$,
$\by^0:= \id$, $\nabla\by^0 := [\be_1,\be_2]$ and $k:=0$. \\
(1) Compute $\theta^{k+1}\in H^1(\o)$ with $\theta^{k+1}|_{\p_D \o} = \theta_D$
and
\[
(\os d_t \theta^{k+1},\vphi) + (\ok \nabla \theta^{k+1},\nabla \vphi) 
+ (\oeta \theta^{k+1},\vphi)_{\p_R \o} = (\oeta \theta_{ext},\vphi)_{\p_R \o}
\]
for all $\vphi \in H^1_D(\o) := \{\psi \in H^1(\o): \psi|_{\p_D \o} = 0\}$. \\
(2) Compute a minimizer $\by^{k+1}\in H^2(\o;\R^3)$ of the functional
\[
I[\by] = \ab{\frac1{12}} \int_\o \omu |D^2\by|^2 - 2 \oa \omu \theta^{k+1} H:I_2
+ 2 \omu (\oa \theta^{k+1})^2
\]
subject to \rhn{the obstacle constraint $\by \in K$,
the} boundary conditions $\by^{k+1} = [\id,0]$ and
$\nabla \by^{k+1} = \nabla [\be_1,\be_2]$ on $\p_D\o$ and the linearized
isometry condition
\[
[\nabla (\by-\by^k)]^\transp \nabla \by^k + [\nabla \by^k]^\transp \nabla (\by-\by^k) =0
\qquad \text{in }\omega.
\]
(3) Increase $k\to k+1$ and continue with~(1). 
\end{algorithm}

Note that the operators $\nabla$ and $D^2$ entail partial
derivatives with respect to the parametric variables
$\bx=(x_1,x_2)\in\omega$ only.
To obtain a practical, semi-implicit version of Algorithm~\ref{alg:alg_abstract}
we use the identities derived in Section~\ref{sec:pdes}.
Given an approximation $\by^k$ we define a tangent space relative to the
isometry constraint and the boundary conditions by 
\[
\cF[\by^k] := \big\{ \bv\in [H^2_D(\o)]^3: 
[\nabla \bv]^\transp \nabla \by^k + [\nabla \by^k]^\transp \nabla \bv =0 \big\},
\]
i.e., vector fields $\bv$ satisfying homogeneous clamped boundary conditions
and the linearized isometry constraint~\eqref{iso}.
Note that the space $H^2_D(\o)$ consists of all functions that vanish
along with their gradients on $\p_D\o$.

\medskip
\rhn{{\bf Obstacle constraint.}
We deal with the constraint $\by\in K$
via a variable splitting and penalization of such splitting in the
$L^2$-norm. We thus  
introduce the auxiliary variable $\bs \approx \by$, penalize the
deviation of $\bs$ from $\by$ by adding the penalty $L^2$-term
$
\ab{\frac{1}{12\veps}} \|\by-\bs\|^2
$
to the energy $I[\by]$, i.e. consider
\[
J[\by,\bs] := I[\by] + \ab{\frac{1}{12\veps}} \|\by-\bs\|^2
\]
and minimize $J[\by,\bs]$ separately over $\by$ and $\bs$
imposing that $\bs\in K$. Note that minimizing with
respect to $\bs$ leads to the $L^2(\o)$-orthogonal projection of $\by$ onto
$K$, which we denote by $\bs = \Pi_K(\by)$, whereas $\by$ is unconstrained.
We assume that the undeformed
plate does not intersect the obstacle.
}

\begin{algorithm}[practical time stepping]\label{alg:alg_practical}
  Let $\tau>0$, set $\theta^0 :=0$, $\by^0 := [\id,0]$,
$\nabla\by^0 := [\be_1,\be_2]$, $\bs^0 := \by^0$, and $k:=0$. \\
(1) Compute $\theta^{k+1}\in H^1(\o)$ with $\theta^{k+1}|_{\p_D\o} = \theta_D$
and 
\[
(\s d_t \theta^{k+1},\vphi) + (\ok \nabla \theta^{k+1},\nabla \vphi) 
+ (\oeta \theta^{k+1},\vphi)_{\p_R \o} = (\oeta \theta_{ext},\vphi)_{\p_R \o}
\]
for all $\vphi \in H^1_D(\o)$. \\
(2) Compute $\bv^{k+1}\in \cF[\by^k]$ such that 
\[\begin{split}
& (\omu D^2[\by^k+\tau \bv^{k+1}],D^2\bw) +  \veps^{-1} (\by^k+\tau \bv^{k+1},\bw) \\
&\hspace*{3cm} =  (\omu \Delta \bw \cdot [\p_1 \by^k \times \p_2 \by^k],\oa \theta^{k+1})
+ \veps^{-1} (\bs^k,\bw)
\end{split}\]
for all $\bw\in \cF[\by^k]$ and set $\by^{k+1} := \by^k + \tau \bv^{k+1}$. Set
\[
\bs^{k+1} := \Pi_K(\by^{k+1}).
\] 
(3) Increase $k\to k+1$ and continue with~(1). 
\end{algorithm}

The precise stability properties of Algorithm~\ref{alg:alg_practical}
appear difficult to identify. In particular, due to the lack of control on a discrete
time-derivative for $(\by^k)$, the inconsistencies related to the explicit
treatment of some terms cannot be controlled directly.
\rhn{If we write $\by=\by^k+\tau\bv$ with $\bv\in\cF[\by^k]$,
then the equation in Step~(2) is the Euler-Lagrange equation
of the energy
\begin{equation}\label{semi-implicit}
\begin{aligned}
J[\by;\by^k,\bs^k,\theta^{k+1}] & = 
\int_\o \Big(\ab{\frac{\omu}{12}} |D^2\by|^2 +\ab{\frac{1}{12\veps}} |\by^{k+1}-\bs^k|^2
\\
& - \ab{\frac{\omu}{6}} \Delta\by
\cdot (\p_1\by^k\times\p_2\by^k)\oa \theta^{k+1}
+ \ab{\frac{\omu}{6}} (\oa\theta^{k+1})^2 \Big) d\bx
\end{aligned}
\end{equation}
without obstacle constraint on $\by$. In contrast, $\bs^{k+1}$
is the $L^2$-projection of $\by^{k+1}$ onto $K$, namely $\bs^{k+1} = \Pi_K(\by^{k+1})$,
which does not involve any solve. The 
decoupling of $\by^{k+1}$ and $\bs^{k+1}$ in Step~(2) is 
motivated by separate convexity properties of $J[\by,\bs]$ in each argument.
A simultaneous minimization of $J[\by,\bs]$
in $\by$ and $\bs$} can be iteratively realized by repeating the two substeps in Step~(2)
where the term $\veps^{-1}(\bs^k,\bw)$ is repeatedly replaced by
$\veps^{-1}(\bs^{k+1},\bw)$ until the subiteration becomes stationary. 

%--------------------------------------------------------------------------------
\subsection{Space discretization}
%--------------------------------------------------------------------------------

Let $\cT_h$ be a partition of the reference domain $\o$ into
quadrilaterals with diameters comparable to $h$. We use the standard
lowest order $H^1$ conforming finite element space $\Q_1$ to discretize the
diffusion equation. We solve approximately the fourth order
nonlinear bending problems with $H^2$ nonconforming discrete Kirchhoff
quadrilaterals \cite{BaBoNo15}. The key idea in their construction is the use of
two $H^1$ conforming finite element spaces to approximate deformations
and deformations gradients together with a
discrete (or reduced) gradient operator $\nabla_h$
that connects the spaces. We
adopt the description of the element from~\cite{BaBoNo15} which
is motivated by the triangular version  considered in 
\cite{Bart13,BaBaHo:80,Brae:07}. 
We let $\Q_r(T)$ and $\mathbb{P}_r(T)$ denote the set of 
polynomials on $T\in\cT_h$ of partial degree $r$ on each variable and
of total degree $r$, respectively. Let $\cN_h$ be the set of
vertices of elements in $\cT_h$ and $\cE_h$ be the set
of edges in $\cT_h$. For every $E\in \cE_h$ we let $\bn_E$ be
a unit normal vector to $E$ and $\bz_E$ be the midpoint of $E$.

\begin{definition}[discrete spaces and operators]
(i) Define the discrete spaces  
\[\begin{split}
\V_h &:= \big\{ w_h \in C(\overline{\o}): 
 \, w_h|_T \in \Q_1(T) ~ \forall T\in \cT_h \big\}, \\
\W_h &:= \big\{ w_h \in C(\overline{\o}): 
 \, w_h|_T \in \Q_3(T) ~ \forall T\in \cT_h, \ 
\nabla w_h \mbox{ continuous in } \cN_h, \\ 
&  \qquad \qquad  \nabla w_h(\bz_E) \cdot \bn_E
= \frac12 \big( \nabla w_h(\bz_{E}^1)+\nabla w_h(\bz_{E}^2)\big) \cdot
\bn_E ~\forall E\in \cE_h, \big\}, \\
\G_h &:= \big\{ \psi_h \in [C(\overline{\o})]^2:  \, \psi_h|_T \in
[\Q_2(T)]^2 ~\forall T\in \cT_h \big\}.
\end{split}\]
(ii) Let $\hcI_h^2: [H^2(\o)]^2\to \G_h$ be the interpolation
operator defined by
\[\begin{aligned}
\hcI_h^2 \bpsi (\bz)  &= \bpsi(\bz) &&  \mbox{ for all $\bz \in \cN_h$}, \\
\hcI_h^2 \bpsi (\bz_E) &= \bpsi(\bz_E) && \mbox{ for all $E\in \cE_h$}, \\
\hcI_h^2 \bpsi (\bz_T) &= \frac14 \sum_{\bz\in \cN_h\cap T} \bpsi(\bz) && \mbox{ for all $T\in \cT_h$}.
\end{aligned}\]
The operator $\hcI_h^2$ is also well-defined for
  discrete vector fields $\psi \in \nabla \W_h$. \\
(iii) Let $\nabla_h : H^3(\o)\to \G_h$  be the discrete gradient
  operator defined by
\[
\nabla_h w := \hcI_h^2 \big[\nabla w\big].
\]
The operator $\nabla_h$ is also well-defined for discrete functions $w \in \W_h$.
\end{definition}

Note that $\nabla\W_h$ is a space of discontinuous vector fields
containing $\Q_2(T)$ for all $T\in\cT_h$, whereas $\nabla_h\W_h$ is
a smaller conforming space.

For an efficient numerical treatment of nonlinearities such as the
projection operator $\Pi_K$ we define the discrete inner product $(\cdot,\cdot)_h$ for 
piecewise continuous functions $\phi,\psi\in \Pi_{T\in\cT_h} [C^0(T)]^\ell$
\[
(\phi,\psi)_h := \sum_{T\in \mathcal T_h}
\frac{|T|}{4} \sum_{\bz \in \cN_h \cap T} \phi|_T(\bz) \cdot \psi|_T(\bz).
\]
To define our numerical scheme we also define
\[
[\p_1^h, \p_2^h] := \nabla_h, \quad D_h^2 := \nabla \nabla_h, \quad
\Delta_h := \trace D_h^2.
\]
The discrete counterpart of the isometry condition
\eqref{isometry} is imposed at the nodes $\cN_h$ of the mesh $\cT_h$, 
which leads to the following discrete linearization:
\[\begin{split}
\cF_h[\by_h] & = \big\{ \bv_h \in [\W_h]^3: \bv_h|_{\p_D \o} = \nabla_h \bv_h|_{\p_D \o} = 0, \\
& \quad [\nabla_h \bv_h(\bz)]^\transp \nabla_h \by_h(\bz)
+ [\nabla_h \by_h(\bz)]^\transp \nabla_h \bv_h(\bz) =0 \quad \forall \bz \in \cN_h \big\}.
\end{split}\]
The following scheme then only requires the solution of linear systems
of equations and the evaluation of low dimensional projections. 

\begin{algorithm}[fully practical scheme]\label{alg:alg_fully}
Let $\tau>0$, set $\theta_h^0 :=0$, $\by_h^0 := [\id,0]$,
$\nabla\by_h^0:=[\be_1,\be_2]$, $\bs_h^0 := \by_h^0$, and $k:=0$. \\
(1) Compute $\theta_h^{k+1}\in \V_h$ with $\theta_h^{k+1}|_{\p_D\o} = \theta_D$
and 
\[
(\os d_t \theta_h^{k+1},\vphi_h) + (\ok \nabla \theta_h^{k+1},\nabla \vphi_h) 
+(\oeta \theta_h^{k+1},\vphi_h)_{\p_R \o} = (\oeta \theta_{ext},\vphi_h)_{\p_R \o}
\]
for all $\vphi_h \in \V_h\cap H^1_D(\o)$. \\
(2) Compute $\bv_h^{k+1}\in \cF_h[\by_h^k]$ such that 
\[\begin{split}
& (\omu D_h^2[\by_h^k+\tau \bv_h^{k+1}],D_h^2\bw_h) + \veps^{-1} (\by_h^k+\tau \bv_h^{k+1},\bw_h) \\
&\hspace*{3cm} =  (\omu \Delta_h \bw_h \cdot [\p_1^h \by_h^k \times \p_2^h \by_h^k],\oa \theta_h^{k+1})_h
+ \veps^{-1} (\bs_h^k,\bw_h)_h
\end{split}\]
for all $\bw_h\in \cF_h[\by_h^k]$ and set $\by_h^{k+1} := \by_h^k + \tau \bv_h^{k+1}$. Set
\[
\bs_h^{k+1}(\bz) := \Pi_K\big(\by_h^{k+1}(\bz)\big)
\]
for all $\bz\in \cN_h$. \\
(3) Increase $k\to k+1$ and continue with~(1). 
\end{algorithm}

The definitions of both $\cF_h[\by_h^k]$ and $\by_h^{k+1}$ imply
that $\by_h^{k+1} = [\id,0], \nabla y_h^{k+1} = [\be_1,\be_2]$
on $\p_D\omega$ for all $k\ge0$.

%%%%%%%%%%%%%%%%%%%%%%%%%%%%%%%%%%%%%%%%%%%%%%%%%%%%%%%%%%%%%%%%%%%%%%%%%%%%%%%%
\section{Numerical experiments}\label{sec:experiments}
%%%%%%%%%%%%%%%%%%%%%%%%%%%%%%%%%%%%%%%%%%%%%%%%%%%%%%%%%%%%%%%%%%%%%%%%%%%%%%%%

We aim at numerically investigating practically relevant scenarios for our
model problem. For this purpose we work with physical units and realistic
material parameters in the following. These apply to polydimethylsiloxane,
polyvinyl-alcohol, polystyrene, soft polyvinylchloride, and soft polyurethane
materials.

%-------------------------------------------------------------------------------
\subsection{Material parameters}
%-------------------------------------------------------------------------------

The material parameters involved in our mathematical model belong
to the following ranges for typical polymer materials:
\[\begin{split}
&\text{thermal conductivity: $\k$  = \SIrange[range-units = single]{0.1}{0.5}{\watt \per \meter \per \celsius}} \\
&\text{material density: $\vrho$ =  \SIrange[range-units = single]{1.0}{2.0e3}{\kilogram \per \cubic\meter}} \\
&\text{specific heat capacity: $c_v$ = \SIrange[range-units = single]{1.0}{2.0e3}{\joule \per \kilogram \per \celsius}} \\
&\text{shear modulus (second Lam\'e parameter): $\mu$ =  \SIrange[range-units = single]{0.2}{2.0e3}{\megapascal}} \\
&\text{first Lam\'e parameter: $\lambda$ = \SIrange[range-units = single]{1.5}{15.0e3}{\megapascal}} \\
&\text{thermal expansion: $\a$ = \SIrange[range-units = single]{\pm
    0.5}{\pm 2.0e-4}{ \per \celsius}} \\
&\text{bilayer thickness: $\d$ = \SIrange[range-units = single]{0.1}{2.0e-3}{\mm}} \\
&\text{heat transfer coefficient (water): $\eta$ = \SI{2.0e6}{\watt \per \square \meter \per \celsius}}
\end{split}\]

The heat capacity is then $\sigma = \vrho c_v$.
% and the diffusivity $\k/\sigma$ 
%belongs to the range \SIrange[range-units = single]{0.1}{0.125}{\square\mm \per \second}.
Except for the thermal expansion coefficient, which has opposite signs in the upper and lower
layer, we use the same material parameters for the two layers. Unless
stated otherwise we use 
\[\begin{split}
&\k=\SI{0.1}{\watt \per \meter \per \celsius},\qquad \quad \lambda = \mu = \SI{1.5e3}{\megapascal}, \\
&\alpha = \pm \SI{0.5e-4}{\per\celsius}, \qquad \eta = \SI{2.0e-3}{\watt \per \square \mm \per \celsius}.
\end{split}\]
We consider different geometries $\o$ whose diameters are on the order
of a few milimeters but with the same thickness
\[
\delta = \SI{1.5e-3}{\mm}.
\]
The resulting effective material parameters for the reduced model are then 
\[\begin{split}
\overline{\k}/\overline{\s} &= \SI{0.1}{\square\mm \per\second}, \quad
\overline{\a} =  \pm \SI{0.1}{\per \mm \per\celsius}, \\
\overline{\mu} &= \SI{2.0e3}{\megapascal}, \quad
\rhn{\overline{\eta}/\overline{\s} = \SI{2.0}{ \mm \per \second}}.
\end{split}\]
Here, we used~\eqref{d:baralph},~\eqref{meanvalues}, 
$\overline{\eta} = (\eta^++\eta^-)/2$, and 
the formula
\[
\overline{\mu} = \mu + \frac{\lambda \mu}{2\mu + \lambda}
\]
from~\cite{FrJaMu02} which results from a more general dimension reduction. 

%-------------------------------------------------------------------------------
\subsection{Discretization parameters}
%-------------------------------------------------------------------------------

\rhn{The choice of meshsize $h$ is dictated by space resolution, which 
improves with the use of mesh refinement in regions of major activity,
typically near the hinges. The choice of time step $\tau$ and penalty
parameters $\veps$ is trickier to yield to realistic physical rather
than numerical effects. To see why, we present now a
nondimensional analysis. We choose characteristic length $\ell$,
temperature $\theta_0$, bending coefficient $\mu_0$ and time $T$, and
denote the scaled variables
\[
\hbx := \ell^{-1}\bx,
\quad \hby(\hbx) := \ell^{-1}\by(\bx),
\quad \hbs(\hbx) := \ell^{-1}\bs(\bx),
\quad \hth (\hbx) := \theta_0^{-1}\theta(\bx),
\]
and scaled parameters
\[
\hmu := \mu_0^{-1}\omu,
\quad \htau := T^{-1}\tau.
\]
We next rewrite the funcional in \eqref{semi-implicit} in terms of the new variables:
\begin{equation*}
\begin{aligned}
J[\by^{k+1} ;\by^k,\bs^k &,\theta^{k+1}] =
\mu_0 \int_{\ho} \Big(\ab{\frac{\hmu}{12}} |\widehat{D}^2\hby^{k+1}|^2
+ \ab{\frac{\ell^4}{12\veps\mu_0}} |\hby^{k+1}-\hbs^k|^2
\\
&
- \ab{\frac{\ell\hmu}{6}} \widehat{\Delta}\hby^{k+1}
\cdot (\widehat{\p}_1\hby^k\times\widehat{\p}_2\hby^k)\oa\theta_0 \hth^{k+1} 
+ \ab{\frac{\ell^2\hmu}{6}} (\oa\theta_0\hth^{k+1})^2 \Big) d\hbx.
\end{aligned}
\end{equation*}
This reveals that the effective penalty parameter is
\[
\heps = \frac{\mu_0}{\ell^4} \veps \ll 1.
\]
Moreover, since the evolution is dictated by the diffusion equation,
we would like the elastic energy minimization to reflect the quasi-stationary
nature of the plate evolution. If the $k$-iterate $\hby^k$ is far from the obstacle,
then $\hbs^k=\Pi_K(\hby^k)=\hby^k$ and we can rewrite the penalty term
as follows
\[
\int_{\ho} \ab{\frac{\ell^4}{12\veps\mu_0}} \big|\hby^{k+1}-\hbs^k\big|^2
= \int_{\ho} \ab{\frac{\ell^4\htau^2}{12\veps\mu_0}} \, \big|d_{\widehat t} \, \hby^{k+1}\big|^2
\]
where $d_{\widehat t} \, \hby^{k+1} = \htau^{-1}\big(\hby^{k+1}-\hby^k\big)$
is the discrete time derivative.
Consequently, to capture the expected physical behavior we
impose the condition
\begin{equation}\label{tau-epsilon}
\frac{\ell^4}{T^2\mu_0} \frac{\tau^2}{\veps} \ll 1
\quad\Rightarrow\quad
\tau^2 \ll \frac{T^2\mu_0}{\ell^4} \veps \ll T^2.
\end{equation}
In other words, the two limits $\tau\to0$ and $\veps\to0$ do not
commute and the former has to take place before the latter for model consistency.
}

%-------------------------------------------------------------------------------
\subsection{Bilayer micro-scale valves and switching devices}\label{subsec:switch}
%-------------------------------------------------------------------------------

\ab{Adaptive bilayer materials are promising for micro-valves for fluid flow management systems, where flapping or flexure mechanisms are used to control the flow. For example, thermally actuated micro-valves made of bilayers \cite{gordon1991thermally} have been proposed where closing and opening a valve aperture is controlled through heating. Another example is using piezoelectric ceramics for micro-valves \cite{li2005development}. As ceramics is brittle, limited deformation can be obtained in closing and opening the valves. Alternatively, polymer bilayers can be actuated to attain relatively large deformations and used in micro-valves.}

\begin{figure}[ht!]
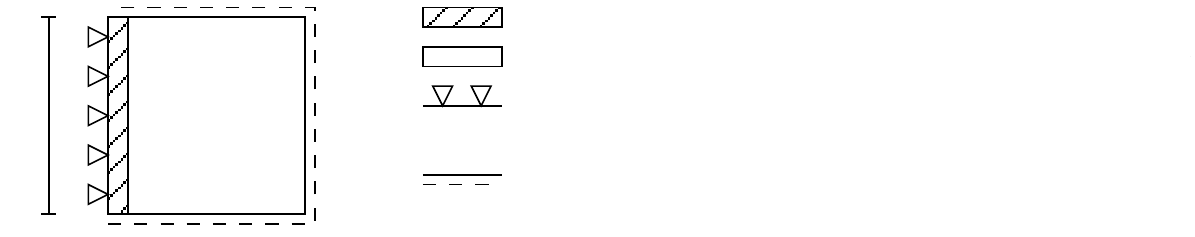 \\[5mm]
\begin{tabular}{lll}
% \hline
\footnotesize $t_0$ & \footnotesize$t_1$  & \footnotesize $t_2$ \\[-1mm]
\includegraphics[width=0.3\textwidth]{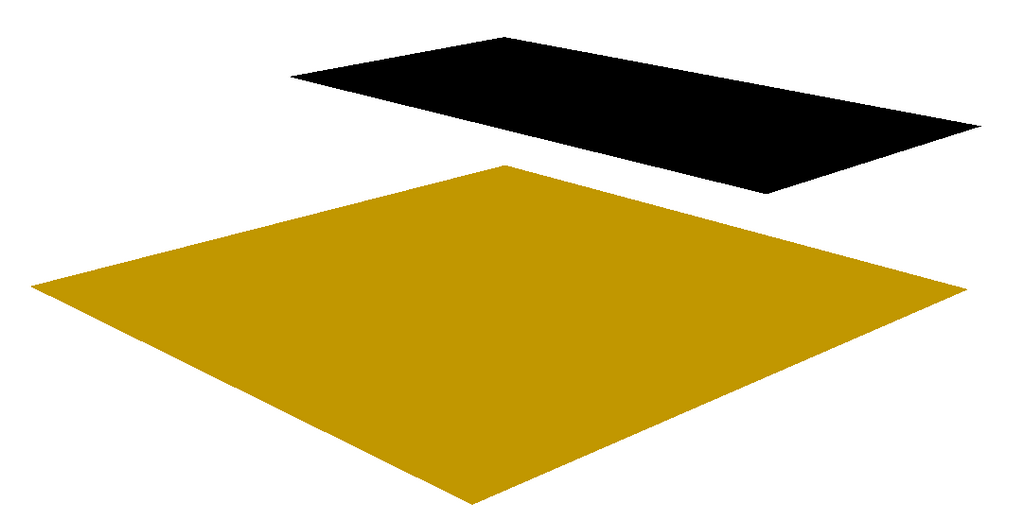}
&
\includegraphics[width=0.3\textwidth]{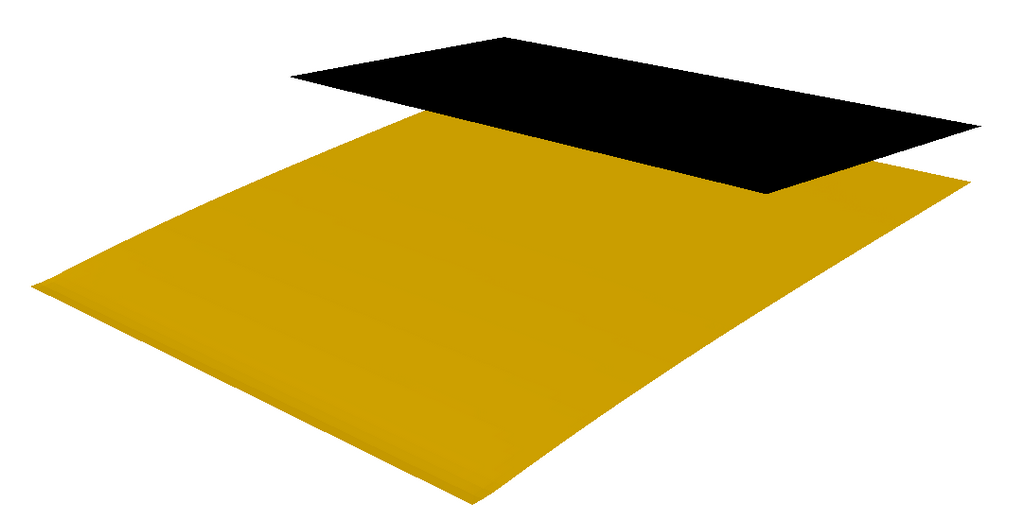}
&
\includegraphics[width=0.3\textwidth]{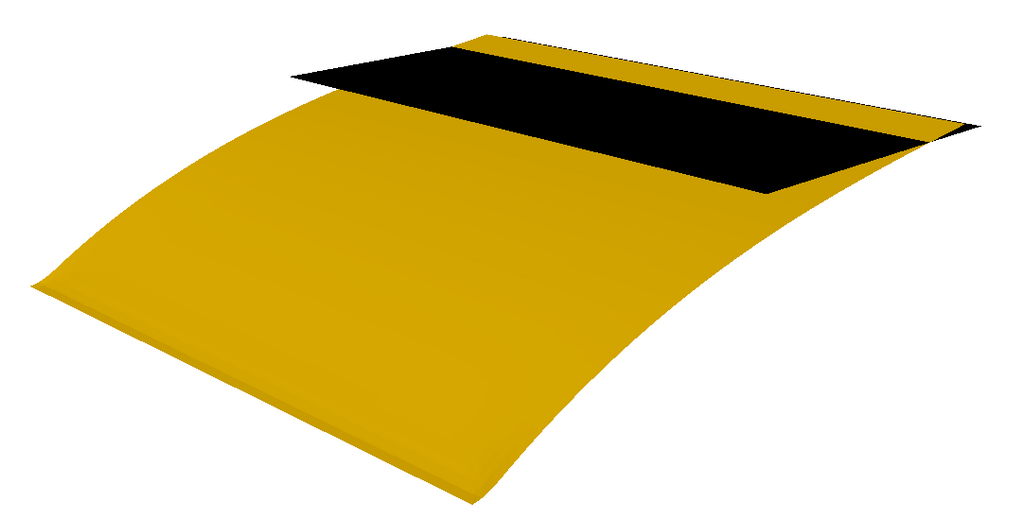}\\
% \includegraphics[width=0.3\textwidth]{Switch6iter400}\\
% \tiny $(a) ~ t=0s$ & \tiny$(b)~t=0.9s$  & \tiny $(c)~t=1.2s$ \\
% \multicolumn{1}{c|}{\includegraphics[width=0.3\textwidth]{Switch6_no_obs}} &
% \includegraphics[width=0.3\textwidth]{Switch6iter3047}&
% \includegraphics[width=0.3\textwidth]{Switch6iter700}\\
% \multicolumn{1}{c|}{ \tiny $(f)~t=+\infty$} & \tiny $(e)~t=+\infty$ & \tiny $(d)~t=2.1s$
\end{tabular}
\begin{center}
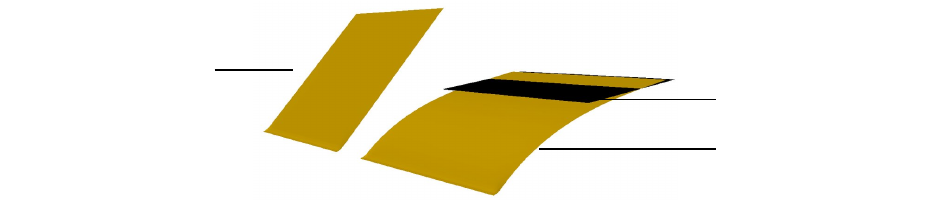 
\end{center}
\caption{\label{fig:bilayer_switch} \label{f:switch_num}
{\em Top:} Bilayer switching device consisting of
a \rhn{flexible single layer and a flexible} bilayer strip that acts as a hinge.  
{\em Middle:} \rhn{Switching device at times 
$t=0.0, 0.9, \SI{2.1}{\second}$. The bilayer
region bends whereas the plate remains flat before hitting the obstacle (a)-(c). 
The device slightly penetrates the obstacle and then bends to accomodate the obstacle. 
{\em Bottom:} Experimental equilibrium states of the bilayer switching device 
without and with the obstacle.}}
\end{figure}

\ab{As proof of concept, we propose in this section a numerical study
of bilayer switching devices consisting of bilayer hinges assembled with flexible plates. Such thermally operated devices have successfully been 
constructed as well}. \rhn{Upon activation the bilayer rotates the plate in order to touch an object.
 We model this by considering the square $\o = (-1,1)^2$ with 
side lengths $\SI{2.0}{\mm}$ which is composed 
of a flexible} single layer and a thin bilayer strip of width $\pi/40$, 
cf.~Figure~\ref{fig:bilayer_switch}. The edge $\p_D \o= \{-1\} \times [-1,1]$ 
of the bilayer strip which is not connected to the plate is assumed
to be clamped, i.e., the deformation $\by$ is constrained to satisfy the boundary 
condition
\begin{equation}\label{clampled}
\by|_{\p_D\o} = [\id, 0]^\transp, \quad
\nabla \by|_{\p_D \o} = [\be_1,\be_2].
\end{equation}
The plate is initially at the critical temperature, i.e., we impose the condition 
\[
\theta(\cdot,0) = \SI{0.0}{\celsius}.
\]
The bending mechanism is triggered by heating the plate from the clamped side
via the prescribed temperature
\[
\rhn{\theta(\cdot,t)|_{\p_D\o} = \min\left(1,\frac t 5\right)~ \SI{100.0}{\celsius}.}
\]
\rhn{The hinge thickness is so small that heat diffuses quite fast. To
prevent a very rapid motion of the plate, which reacts instantaneously
to the hinge bending, the factor $ \min\left(1,\frac t 5\right) $
ensures a smooth and slow transition from the initial temperature
$\theta = \SI{0}{\celsius}$ to the desired one of $\SI{100}{\celsius}$.}
We assume that the plate is thermally insulated and that the deformation is
free on the remaining part $\p_R\omega:=\p\o\setminus \p_D \o$ of the
boundary, i.e., we use
a vanishing effective heat transfer coefficient $\overline{\eta}$.
\rhn{The plate rigidity is uniform throughout the plate and characterized by a bending coefficient $\overline{\mu}=\SI{2.0e3}{\megapascal}$}.
A flat obstacle that models the contact object is modeled by the 
set
\[
\rhn{K = \big\{ \by=(y_i)_{i=1}^3 \in H^1(\o;\R^3): y_3 \le \SI{0.5}{\mm} \big\},}
\]
and the projection $\Pi_K$ onto $K$ is computed nodewise
for $\by_h=(y_i)_{i=1}^3$:
\[
\Pi_K(\by_h)(\bz) =
\big(y_1(\bz),y_2(\bz),\min\{y_3(\bz),0.5\}\big)
\quad\forall \, \bz\in\cN_h.
\]
For our numerical experiments we use a partition of $\o$ that results
from~6 uniform quad refinements of $\o$.
\rhn{We use uniform time steps $\tau = \SI{3.0e-3}{\second}$ 
and penalization parameter $\varepsilon = \SI{4e-6}{mm^4  \megapascal^{-1}}$.
Since $\mu_0=2\times 10^3\megapascal$ and $T=10 s$,
this choice is consistent with \eqref{tau-epsilon}, namely
$\tau^2 \ll 0.5\times 10^{-1} \, s^2$.
Figure~\ref{f:switch_num} shows the evolution from the initial configuration.\\}

\modif{
To assess the intricate influence of the penalty parameter $\varepsilon$, we consider the same setting except for values of $\varepsilon =  4\times 10^{-j}\,\SI{}{mm^4  \megapascal^{-1}}$ for $j=4,...,9$.
Cuts along the plane $\{ x_2 = 0 \}$ (perpendicular to the obstacle and the clamped side) of the plates stationary deformations are displayed in Figure~\ref{fig:switch_cuts}. 
Plates are considered in a stationary state whenever
\begin{equation}\label{e:switch_stat}
\| \by_h^{k+1} - \by_h^k \|_{L_2(\omega)} + \| \nabla \nabla_h (\by_h^{k+1} - \by_h^k) \|_{L_2(\omega)} \leq  10^{-5}.
\end{equation}
On the one hand, the smallest value $\varepsilon$ corresponding to $j=9$ results in a plate which does not cross the obstacle but is influenced by the presence of the obstacle already when still far from the plate.
On the other hand the largest value of $\varepsilon$ corresponding to $j=6$ results in a effective deformations without obstacle. 
Values of $\varepsilon$ for parameters $j=7,8$  yield a deformation crossing the planar obstacle by no more that the finite element meshsize.

\begin{figure}[ht!]
\includegraphics[width=0.5\textwidth]{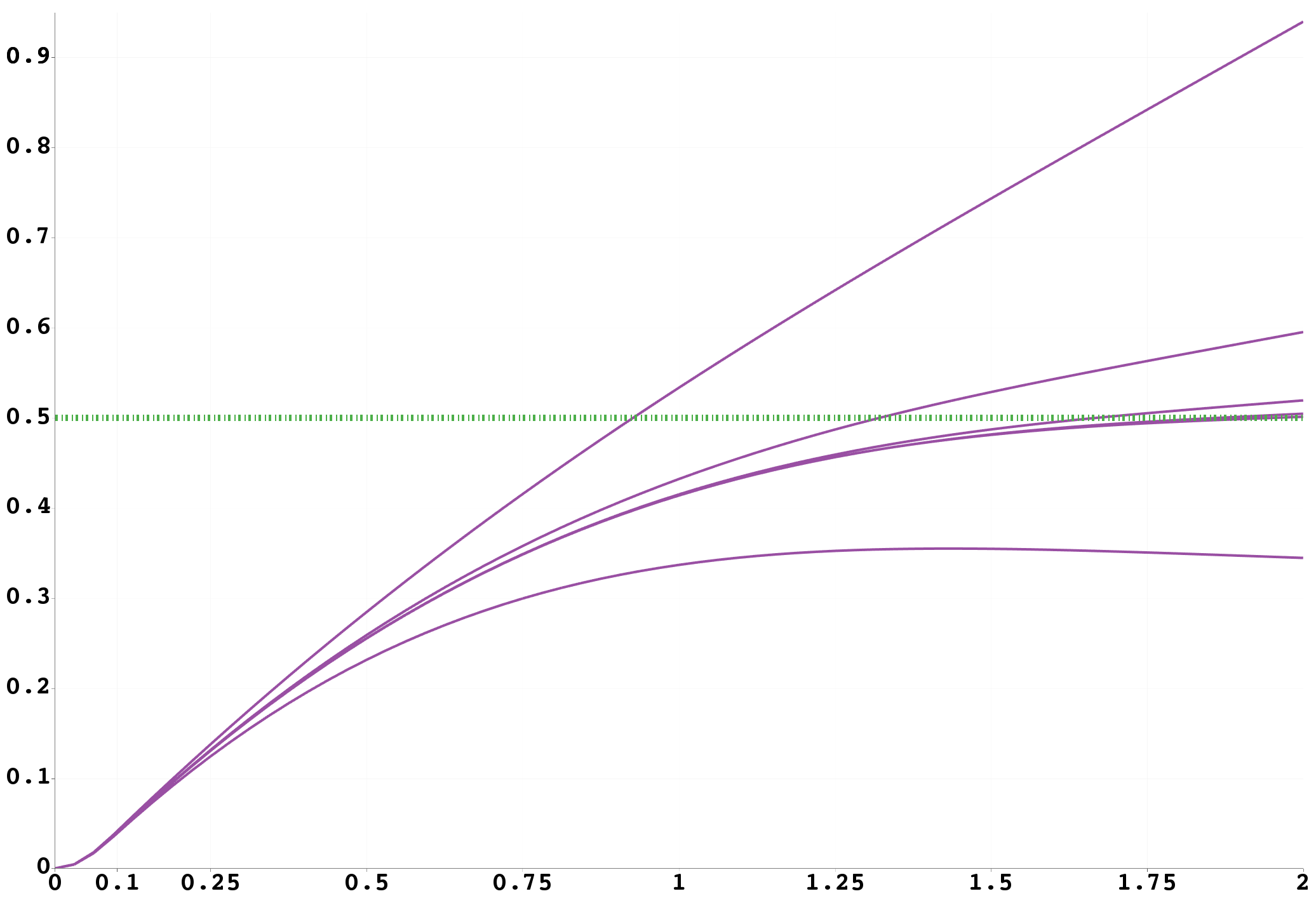} 
\caption{\label{fig:switch_cuts}
\modif{Cuts of stationary deformations along the plane $\{ x_2 = 0 \}$ of the bilayer switching device for (top to bottom curves) $\varepsilon = 4 \times 10^{-j} \SI{}{mm^4  \megapascal^{-1}}$, $j=4,...,9$. Note that the curves corresponding to $j=8$ and $j=7$ are barely distinguishable and do not cross the obstacle $\{x_3 =0.5\}$ more than the finite element meshsize $(=1/64)$. The horizontal line corresponds to the position of the obstacle.}
}
\end{figure}

It is worth mentioning that except the case $j=9$, all the stationary states (defined according to \eqref{e:switch_stat}) are reached at $T=\SI{136.67}{s}$ (44689 time iterations) while $j=9$ required $T=\SI{182.739}{s}$ (60913 time iterations). However, the dynamics for smaller $\varepsilon$ is influenced by the obstacle at early stages as illustrated in Figure~\ref{fig:switch_dyn}.
This is consistent with \eqref{tau-epsilon}. 

\begin{figure}[ht!]
\includegraphics[width=0.3\textwidth]{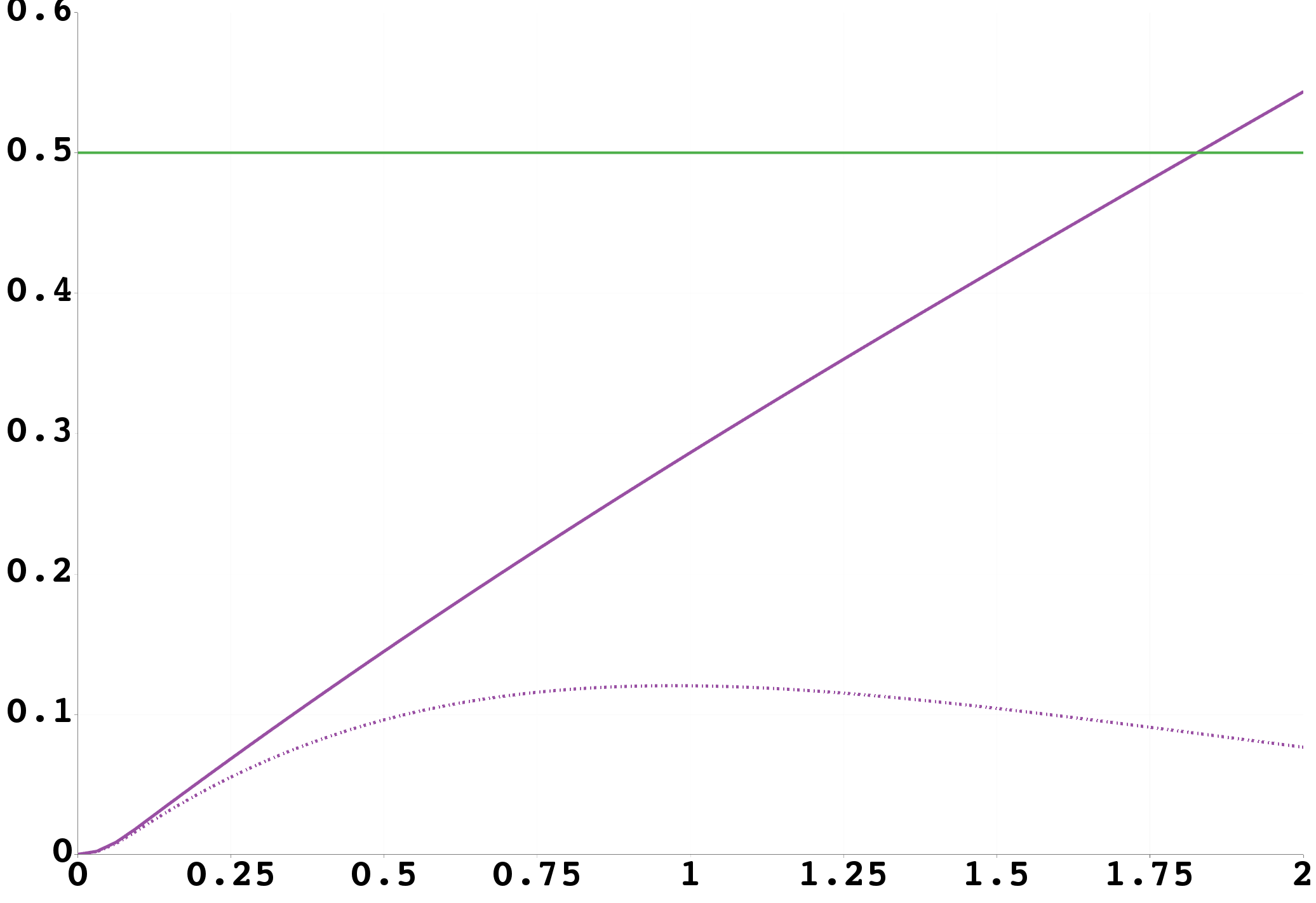} 
\includegraphics[width=0.3\textwidth]{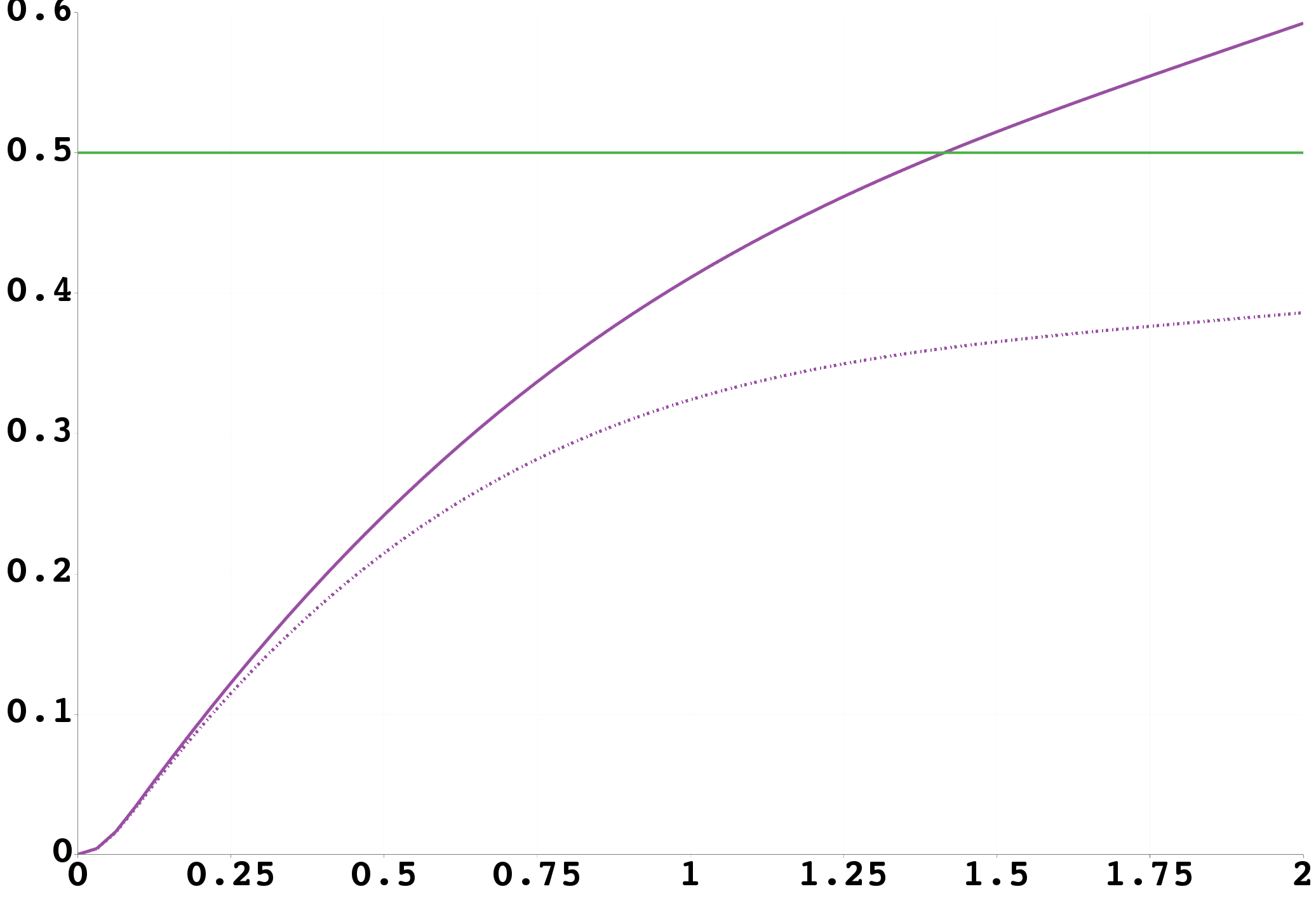} 
\includegraphics[width=0.3\textwidth]{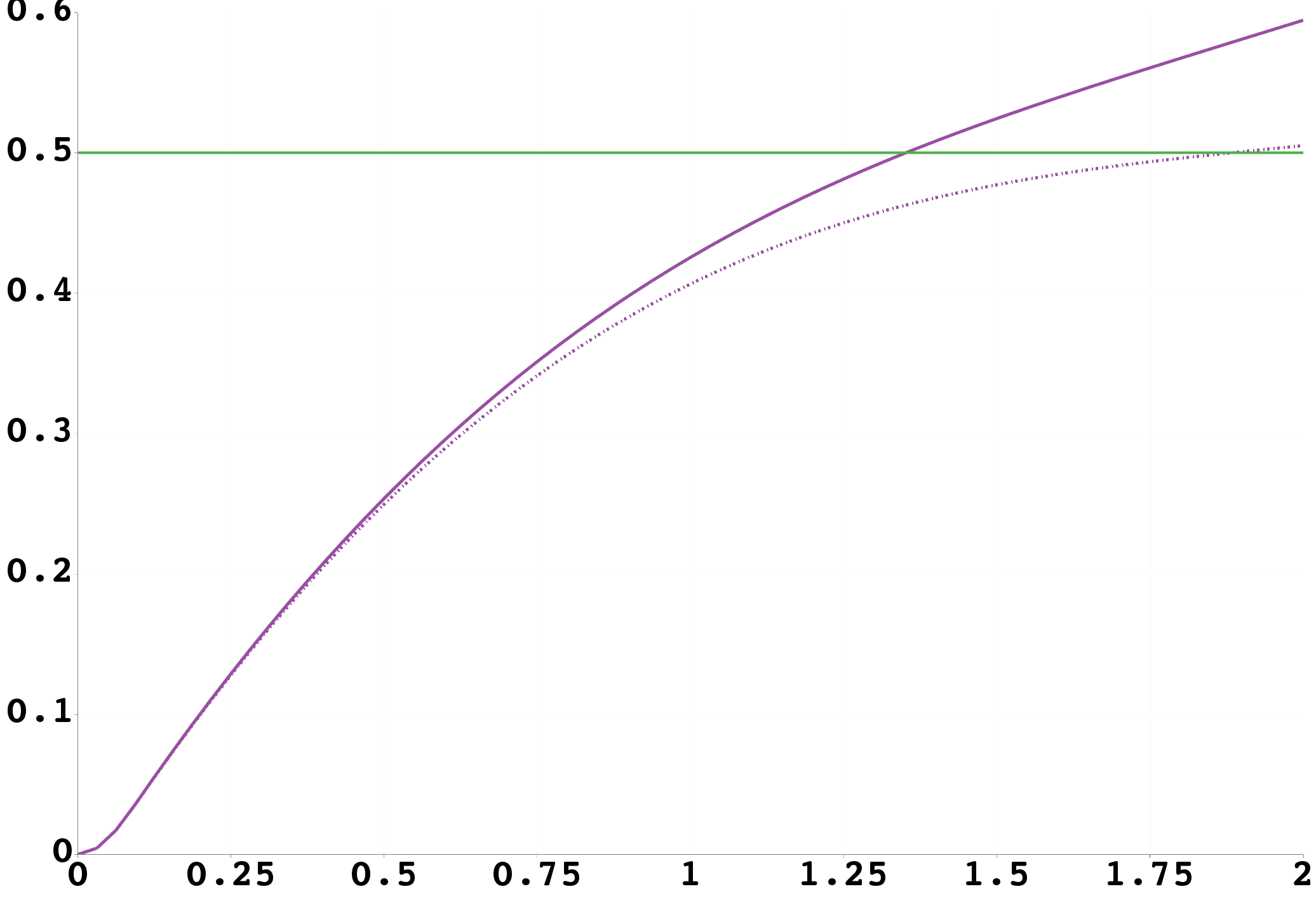} 
\caption{\label{fig:switch_dyn}
\modif{Cuts of deformations along the plane $\{ x_2 = 0 \}$ of the bilayer switching device for (top curve) $\varepsilon = \SI{4e-5}{mm^4  \megapascal^{-1}}$ and (bottom curve) $\varepsilon = \SI{4e-7}{mm^4 \megapascal^{-1}}$ and at times (left) $t=\SI{2.4}{s}$, (middle) $t=\SI{4.8}{s}$ and (right) $t=\SI{7.2}{s}$.}
}
\end{figure}

}

%--------------------------------------------------------------------------------
\subsection{Dog-ear formation}
%--------------------------------------------------------------------------------

\begin{figure}[ht!]
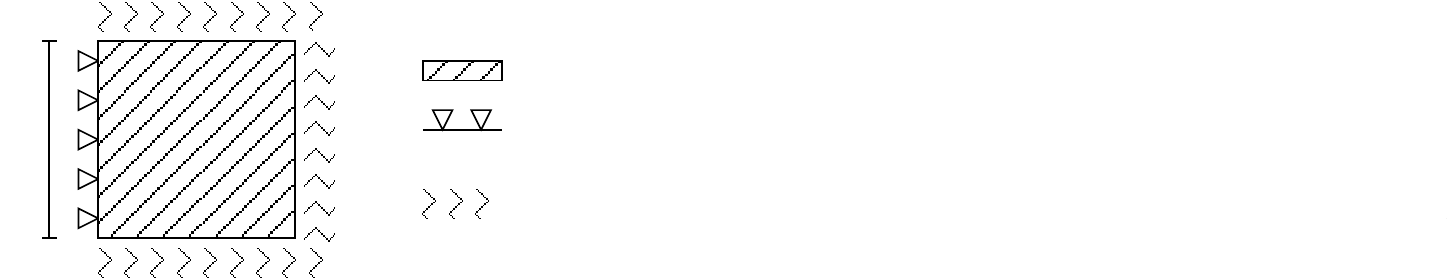 \\[5mm]
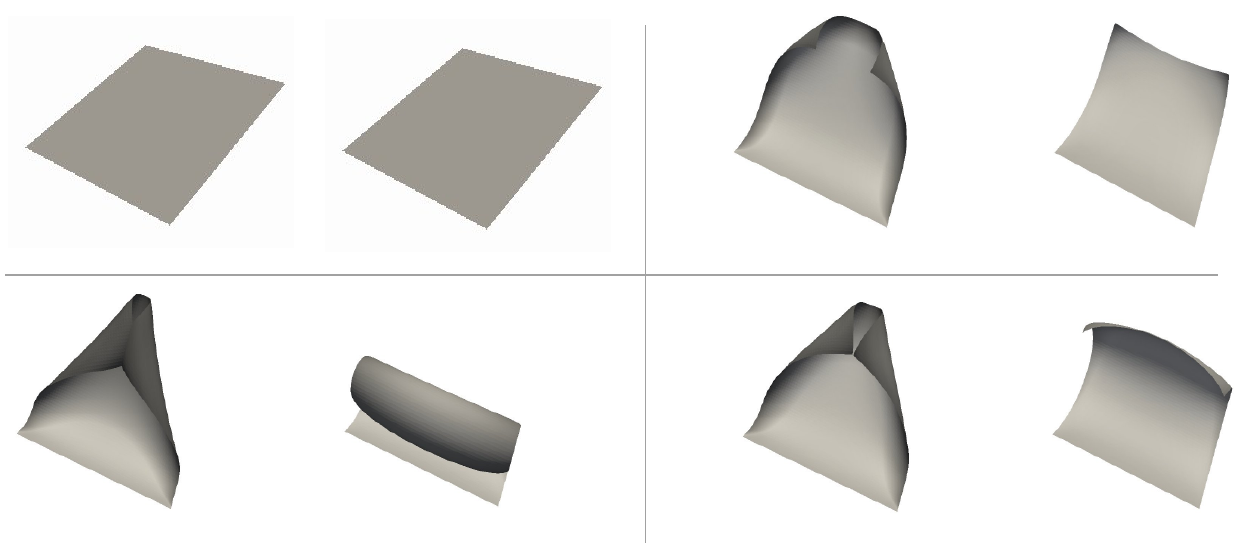
\caption{\label{fig:bilayer_dogear} \label{f:dog_ear}
{\em Top:} Bilayer plate clamped along one side $\p_D\o$
and heated via temperature difference with environment on $\p_R\o=\p\o\setminus\p_D\o$.
{\em Bottom:} Bilayers heated via Robin condition on $\p_R\o$ 
with small~(a) and large~(b) diffusivity at different times.
\rhn{The gray scale represents the temperature within the plate \modif{(brightest = \SI{0}{\celsius}; darkest = \SI{50}{\celsius})}.}
Small diffusivity leads to heat concentrations near $\p_R\o$ and 
strong localized bending behavior which resemble dog-ear formation.}
\end{figure}

To explore possible failure of controlled production of microtubes
reported in~\cite{SZTDI12,Yeetal15,Yeetal16} 
we consider a bilayer square $\o=(-1,1)^2$ of side lengths $\SI{2.0}{\mm}$
that is clamped and thermally insulated on the side
$\p_D\o = \{-1\}\times [-1,1]$, namely $\by_h$ satisfies
  \eqref{clampled} and $\partial_\bn\theta = 0$ on $\p_D\o$.
On the remaining part $\p_R \o = \p\o\setminus  \p_D \o$ the plate is
free and heated via external heat transfer obeying Newton's law of cooling, i.e., 
the normal flux is proportional to temperature difference 
\[
\overline{\k} \nabla \theta \cdot \bn = \overline{\eta}(\theta_{ext} - \theta)
\qquad\text{on } \p_R\o
\]
relative to the ambient temperature $\rhn{\theta_{ext} = \SI{50.0}{\celsius}}$.
The experiment setup is sketched in Figure~\ref{fig:bilayer_dogear}.
For our numerical experiments we consider two different values of 
the effective thermal conductivity $\ok$ and heat capacity $\overline{\sigma}$
so that the resulting diffusivity differs by a factor~10, i.e., 
we consider
\[
\text{(a)} \quad \overline{\kappa}/\overline{\sigma}=\SI{0.1}{\square \mm \per \second}, \qquad
\text{(b)} \quad \overline{\kappa}/\overline{\sigma}=\SI{1.0}{\square \mm \per \second}.
\]
We run Algorithm \ref{alg:alg_fully} (fully discrete scheme) for
a uniform partitition of $\o$ into subsquares with side lengths
$h = 2^{-6} \times \SI{2.0}{\mm}$ and uniform time step $\tau = \SI{5.0e-3}{\second}$.
Figure~\ref{f:dog_ear} displays snapshots of the approximate evolution 
for settings~(a) and~(b) for the respective times
\[
t = 0.0,\, 1.0, \, 2.5, \, 16.0 \times  \overline{\k}/\overline{\sigma}.
\]

We observe
that the smaller diffusivity case (a) leads to accumulation of heat
and higher temperatures in the vicinity the boundary $\p_R\o$ and in particular at
the two corners belonging only to $\p_R\o$. This results in a
stronger localized bending behavior for (a) and the formation of dog-ears which 
then prevent a controlled rolling up into a cylindrical shape. For larger
diffusivity (b), heat distributes faster within the plate and the bending behavior
is less localized.

%---------------------------------------------------------------------------------
\subsection{Self-assembling box}\label{S:self-assembling-box}
%---------------------------------------------------------------------------------

\begin{figure}[ht!]
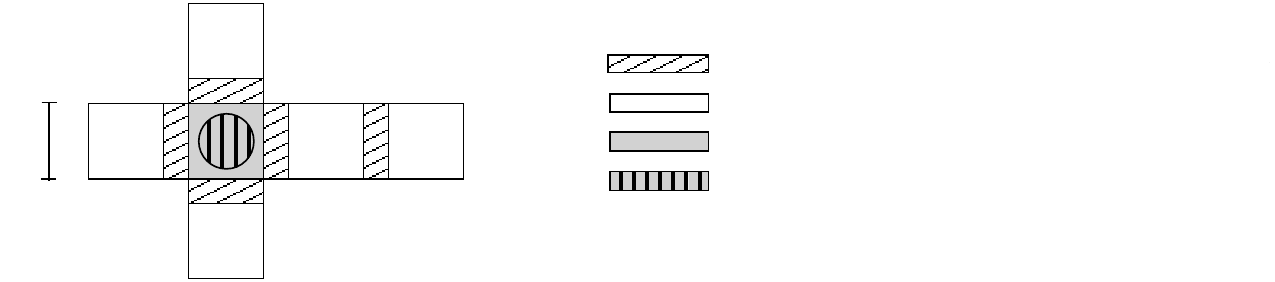 \\[5mm]
\begin{tabular}{ccc}
\footnotesize $t_0$ & \footnotesize $t_1$ & \footnotesize $t_2$ \\
 \includegraphics[width=0.3\textwidth]{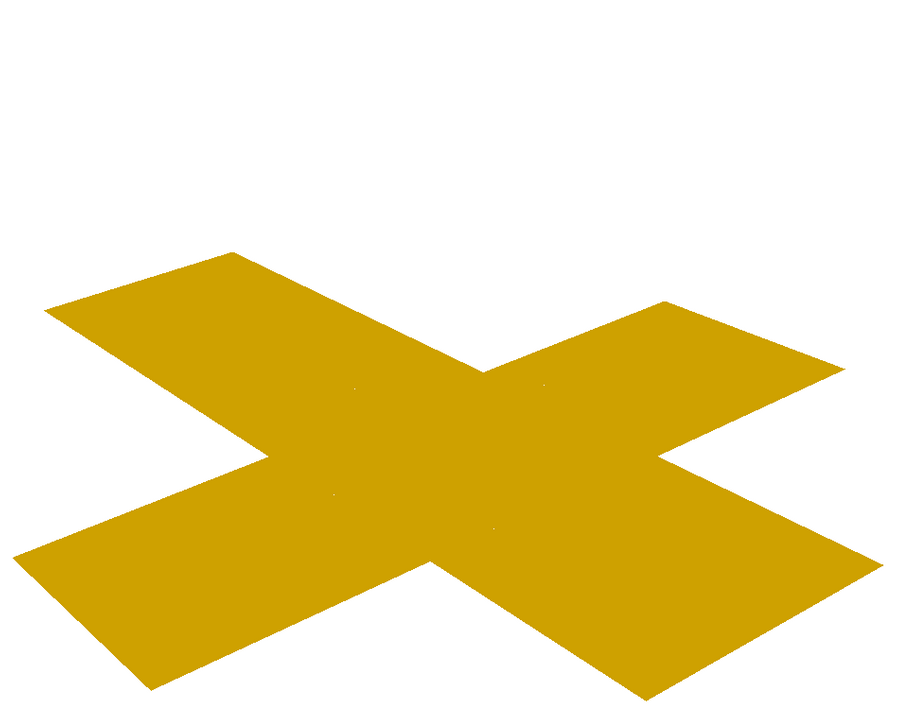} &
 \includegraphics[width=0.3\textwidth]{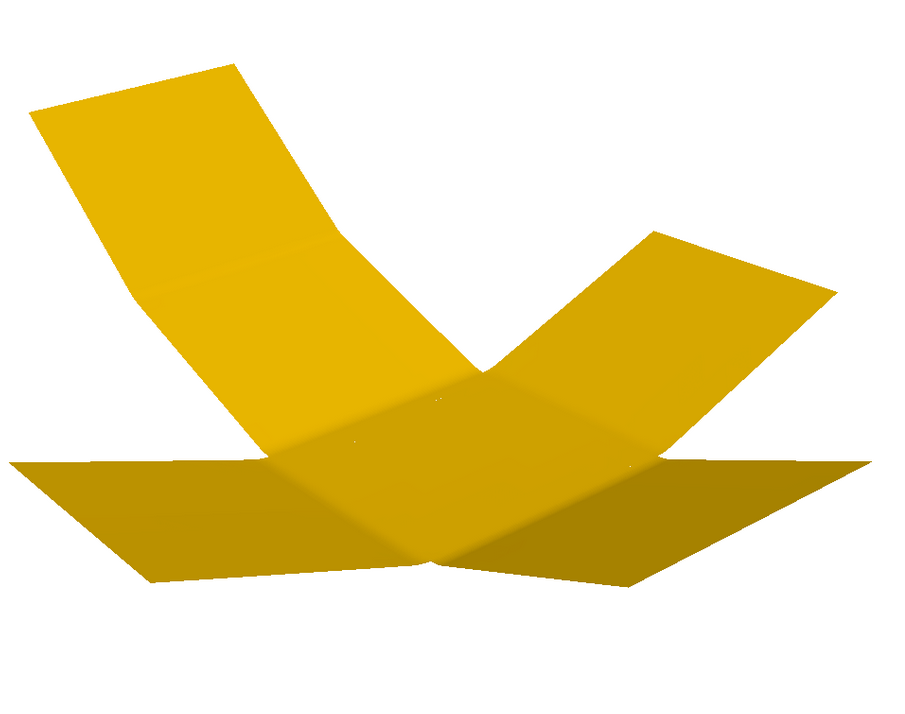}&
 \includegraphics[width=0.3\textwidth]{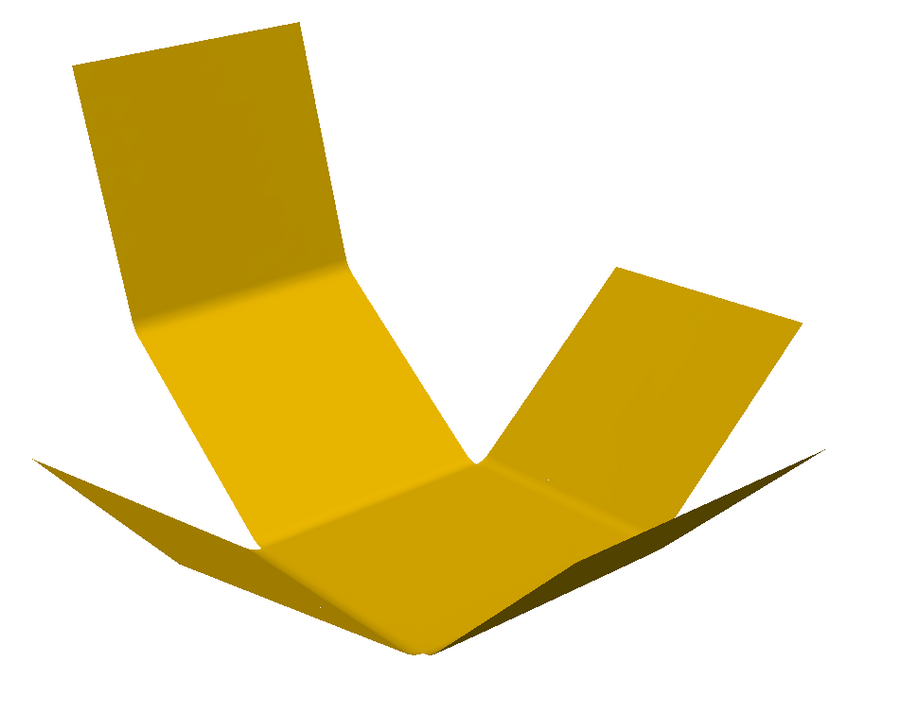} \\
% \tiny $t=0s$ & \tiny $t=5s$ & \tiny $t=10s$ \\
 \includegraphics[width=0.3\textwidth]{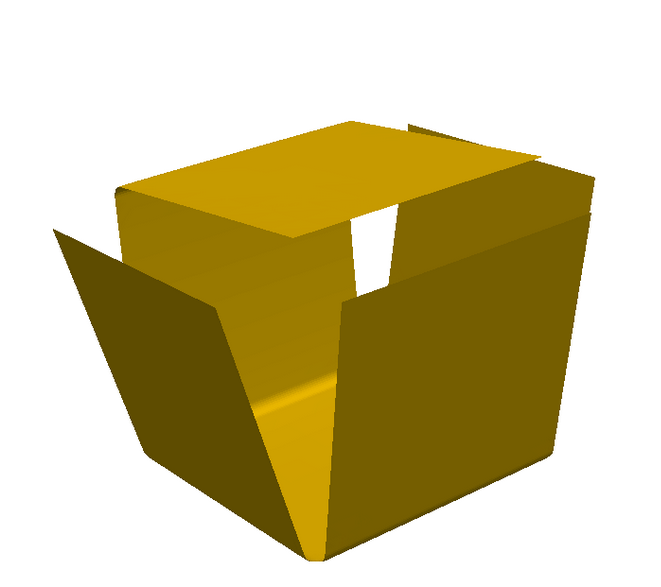}&
 \includegraphics[width=0.3\textwidth]{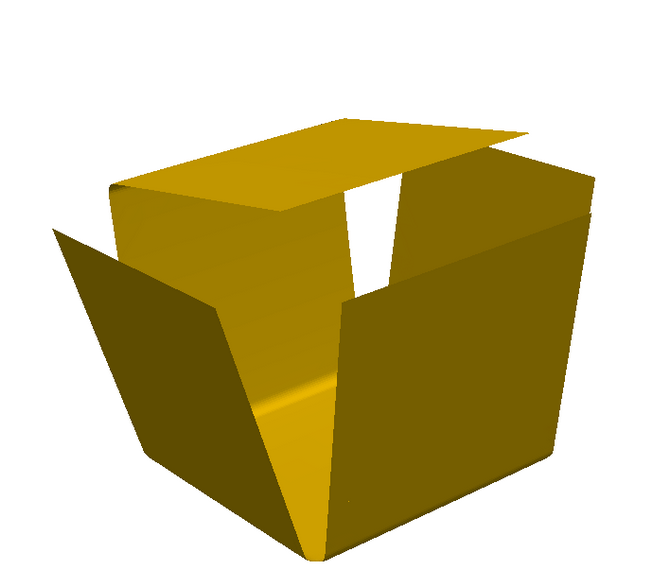}&
 \includegraphics[width=0.3\textwidth]{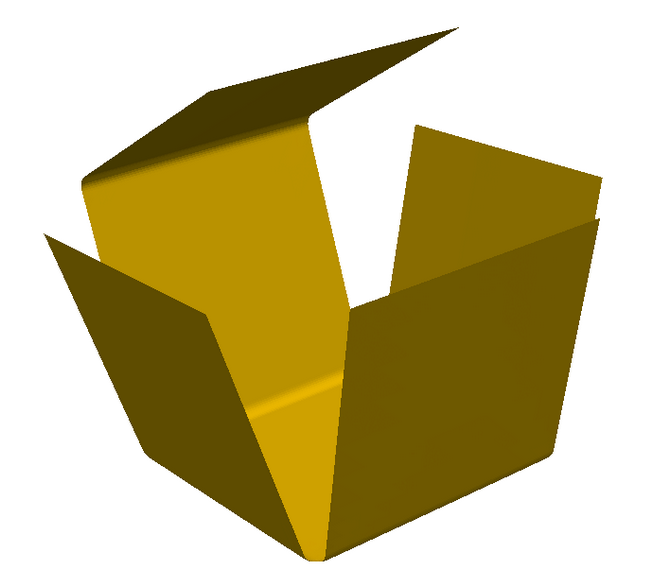}\\
%  \tiny $t=30s$ & \tiny $t=25s$ & \tiny $t=19s$ 
\footnotesize $t_5$ & \footnotesize $t_4$ & \footnotesize $t_3$ 
 \end{tabular}
\caption{\label{fig:self_ass_box} \label{f:box}
{\em Top:} Construction of a self-assembling
box with~6 rigid plates connected by bilayer hinges.
{\em Bottom:} Snapshots after $t=0.0, 5.0, 10.0, 19.0,
25.0, \SI{30.0}{\second}$ of a self-assembling composite box
at different times.
The temperature distribution is nearly uniform throughout the
evolution.} 
\end{figure}

We consider an arrangement of rigid plates and bilayer hinges
that realizes a  
self-assembling and self-opening box actuated by temperature. 
We refer the reader to~\cite{SuShMi94,C6MH00195E} for related
experiments. The relevant practical applications on this type of folding systems are on deployable structures such as drug delivery systems, deployable shelters with prefabricated walls.
As in Subsection~\ref{subsec:switch}, the plate has
$\oa=0$ and $\overline{\mu}$ that is $20$ larger than that of the hinge,
\rhn{so as to mimic a rigid plate insensitive to temperature.}
The geometry is sketched in Figure~\ref{fig:self_ass_box}.

The square rigid plates have side lengths $\SI{1.0}{\mm}$ and
the hinges have widths $\SI{\pi/48}{\mm}$. 
The central square, which is connected to four other squares, 
is assumed to be permanently attached to a substrate, and 
a constant heat source of $\SI{75.0}{\celsius \per \second}$ is applied
on a circle of radius $\SI{0.25}{\mm}$ centered on this 
square during a period of $\SI{19.0}{\second}$.
Note that in that case, the heat equation reads
\[
\partial_t \theta - 10 \Delta \theta = \left\lbrace
\begin{array}{ll}
75.0 &\qquad \text{in }\omega \times (0,19), \\
0.0  & \qquad \text{in }\omega \times (19,\infty).
\end{array}\right.
\]
When temperature increases the bilayer hinges start bending 
while the rigid plates remain undeformed, cf.~Figure~\ref{f:box}.
For sufficient heating a box is formed which can be unfolded via 
external cooling. For our simulation we chose a
uniform time step $\tau=\SI{0.5}{\second}$
and a partition of $\o$ into $3$ uniform quad 
refinements of each square plate, followed by $4$   
quad refinements of elements in the hinges and additional 
local refinements to have at most one hanging node per side.

%---------------------------------------------------------------------------------
\subsection{Deployable airfoil}
\begin{figure}[ht!]
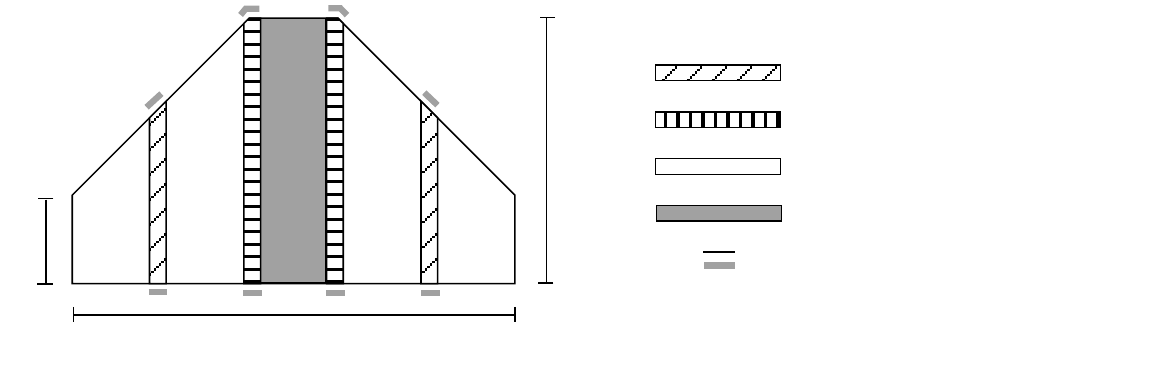 \\[2mm]
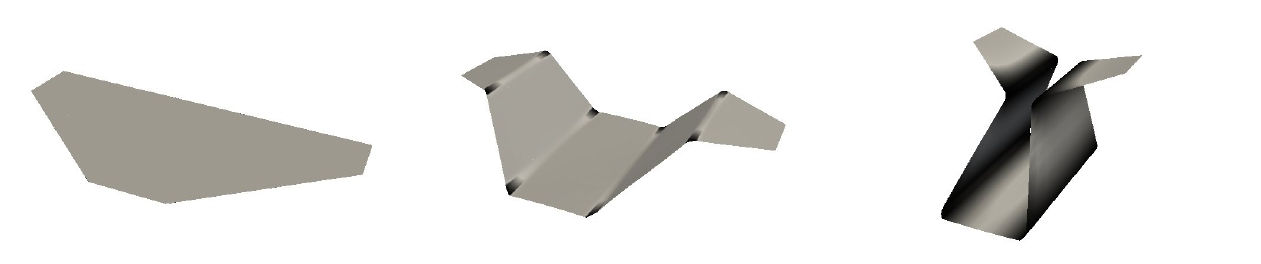
\caption{\label{f:airfoil}\label{fig:airfoil}
{\em Top:} Arrangement of bilayers and rigid plates for the construction
of a deployable airfoil.
{\em Bottom:} 
Snapshots of the folding process of the deployable airfoil
at times $t=$ 0.0, 2.5, 500\SI{}{\second}. The gray scale represents
the temperature within the device \modif{(brightest = \SI{0}{\celsius}; darkest = \SI{60}{\celsius})}.}
\end{figure}

%---------------------------------------------------------------------------------
We consider another arrangement of bilayer hinges and rigid plates
shown in Figure~\ref{fig:airfoil} that models an initially flat
device that can 
fold into an airfoil via external heating. Radical shape changes/folding of wing structures during flight can lead to high maneuverability of the plane and an efficient cruising that is adjustable to various flight envelopes. One example of morphing wing technology, i.e. z-shape folding, has been proposed by the Defense Advanced Research Projects Agency (DARPA), \cite{love2007demonstration}.

The physical parameters are the same as in Subsection
\ref{S:self-assembling-box} \rhn{except for $\overline{\kappa}/\overline{\sigma}=\SI{0.1}{\mm \square \per \second}$} and the middle plate is fixed.
As indicated in the sketch of Figure~\ref{f:airfoil},
the bilayer hinges bend into different directions, which is modeled
by different signs of the effective expansion coefficients,
i.e. $\overline \alpha = \pm 0.3$, or equivalently upon inverting
the bilayer.
The heat diffusion process is initiated by prescribing the 
temperature $\theta_D= \SI{60.0}{\celsius}$ at the Dirichlet boundary parts
$\p_D\o$ where the hinges meet $\p\o$, and vanishing Neumann condition
on the rest $\p_R\o=\p\o\setminus\p_D\o$ (i.e. $\overline{\eta}=0$).
A simulation of the thermally driven folding process is 
illustrated in Figure~\ref{f:airfoil}.

%------------------------------------------------------------------------------
\subsection{Particle encapsulation}
%------------------------------------------------------------------------------

Thermally controlled bilayers for
transporting particles at microscales have been tested experimentally
in~\cite{StPuIo11}. Once deployed in a desirable place,
enclosed particles may be released via external cooling of the device. 
Similar mechanisms can be triggered by changing the pH concentration
of a surrounding liquid. 
This may find exciting and important medical applications
in targeted drug delivery.

\begin{figure}[ht!]
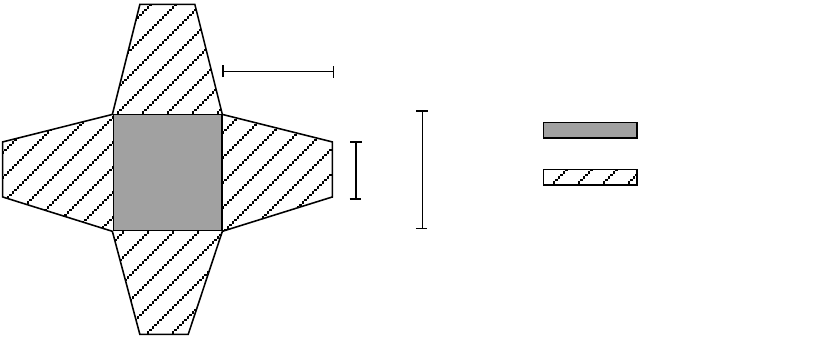 \\[1mm]
\begin{tabular}{cc}
\includegraphics[width=0.3\textwidth]{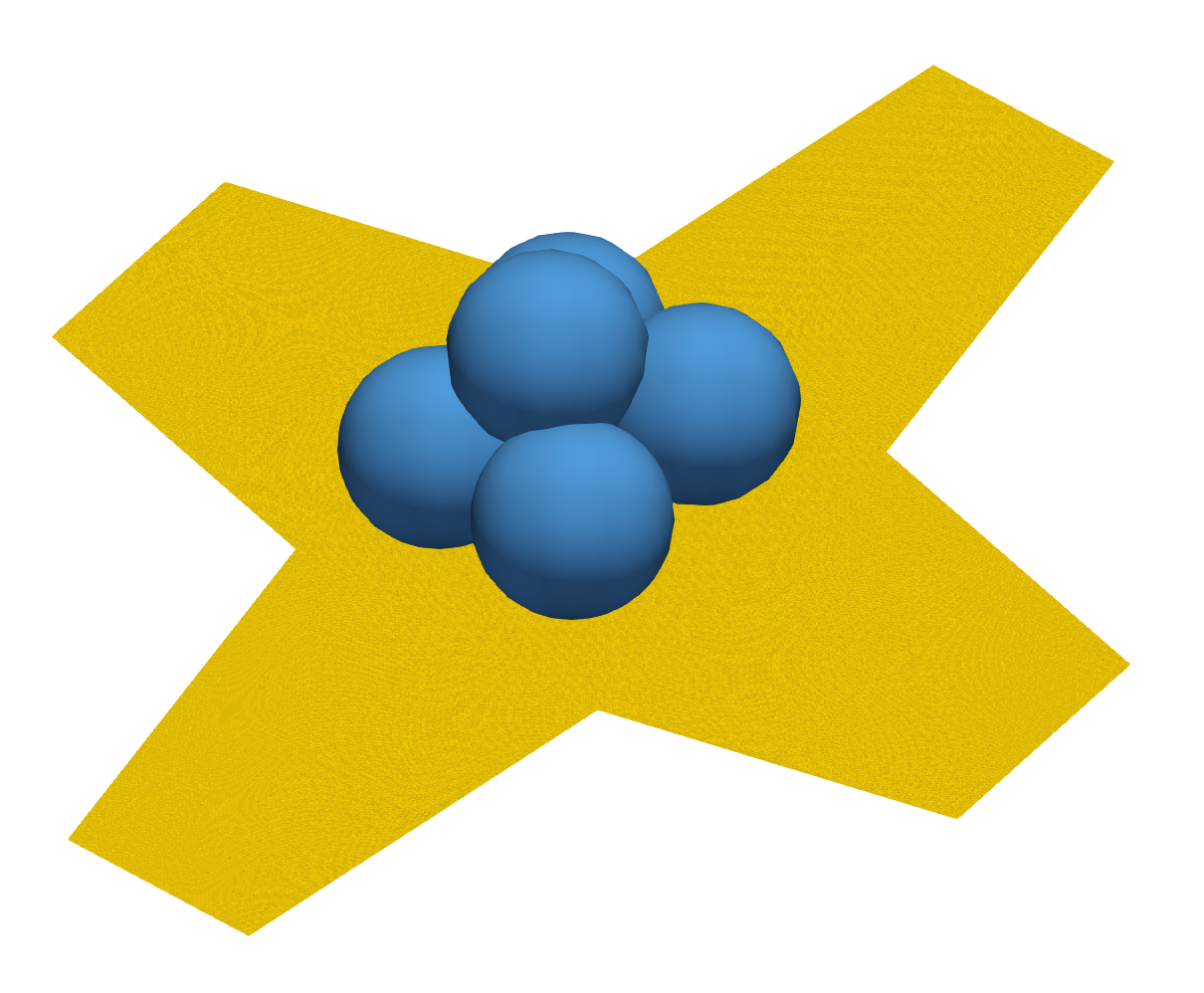}&
\includegraphics[width=0.3\textwidth]{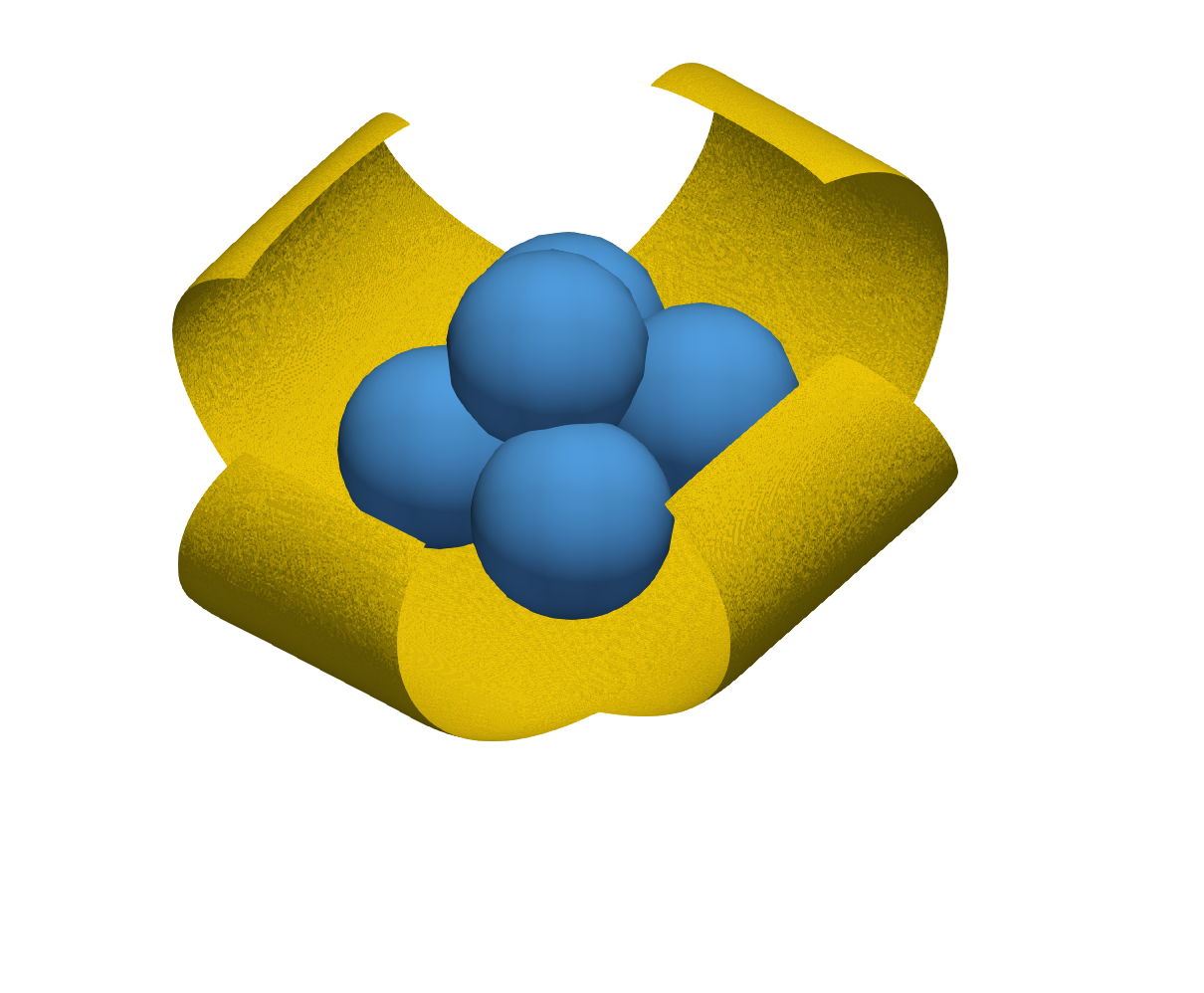}\\
\includegraphics[width=0.3\textwidth]{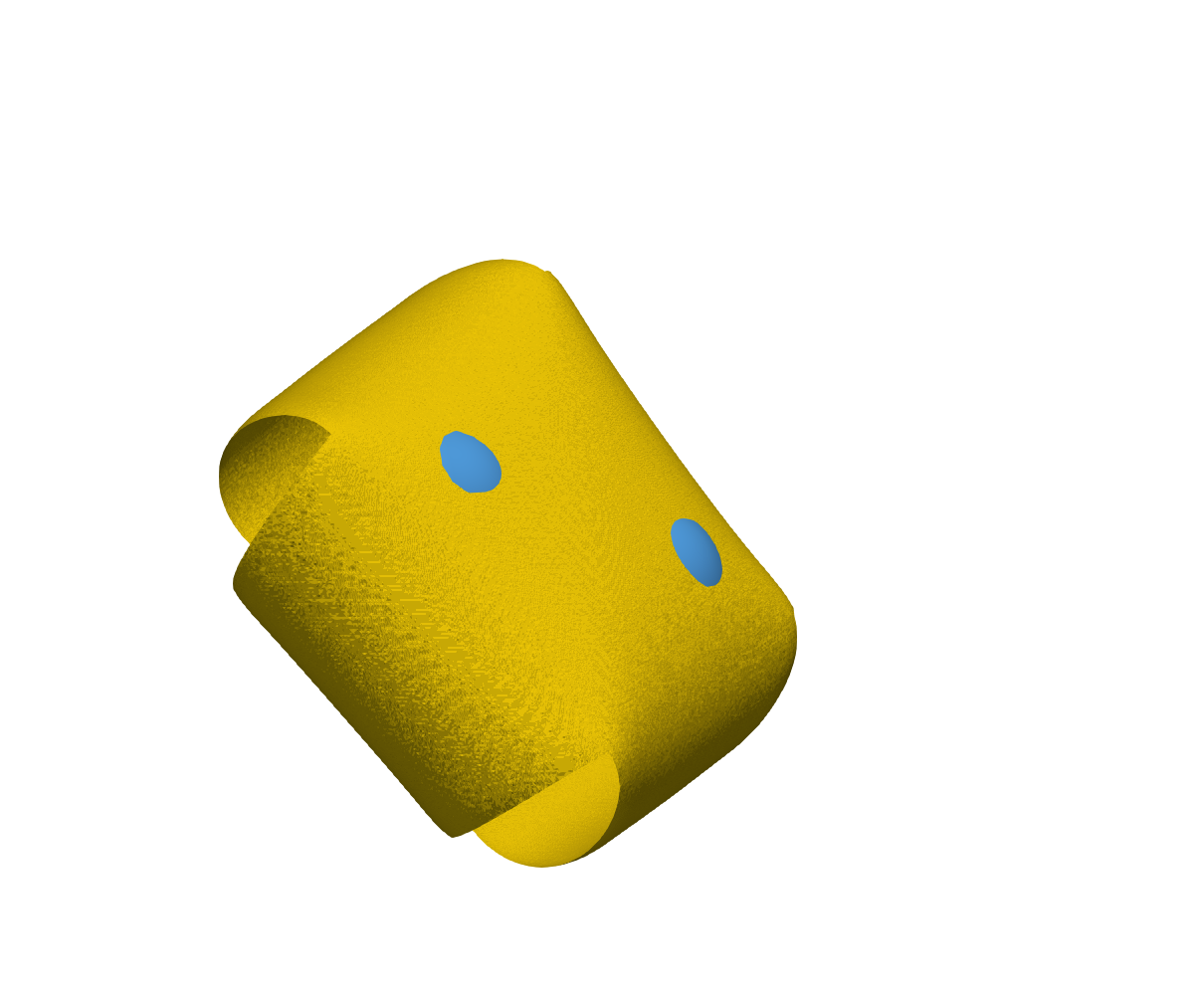}&
\includegraphics[width=0.3\textwidth]{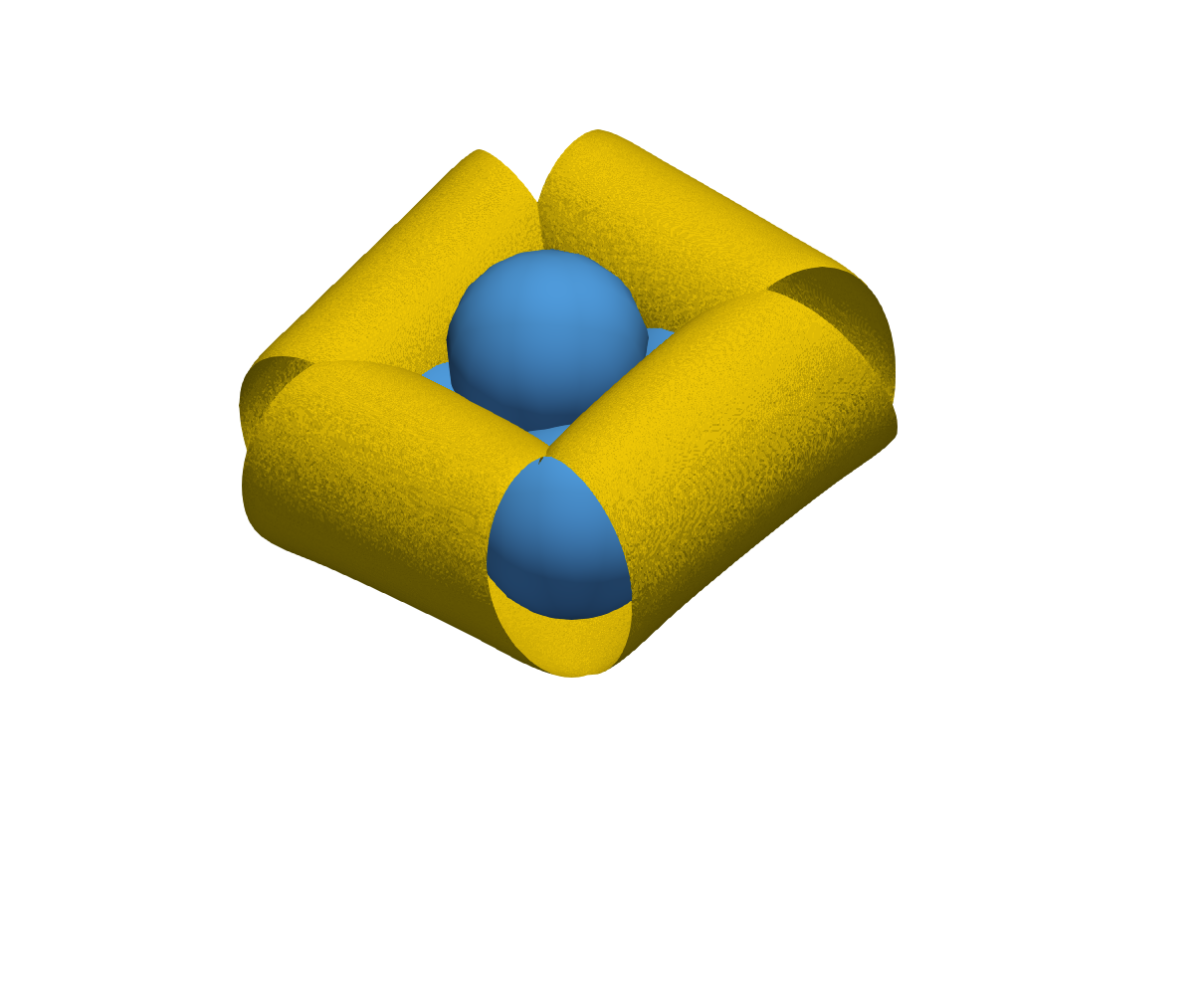}
\end{tabular}
\begin{picture}(0,0)
\put(-4000,1400){$t_0$}
\put(-1800,1400){$t_1$}
\put(-1800,-200){$t_2$}
\put(-4000,-200){$t_2$}
\end{picture}
\caption{\label{fig:microcap} \label{f:microcap_sim}
{\em Top:} Geometry of a bilayer plate that serves as
a microcapsule and can be closed and opened via external heating 
and cooling. The dashed lines indicate the coarsest partition.
{\em Bottom:} Folding of a micro-capsule enclosing~5 non-penetrable particles 
at times $t_0=\SI{0.0}{\second}$, $t_1=\SI{0.5}{\second}$ 
\ab{and} $t_2=\SI{1.0}{\second}$.
The bottom \modif{left} plot provides a different view at $t_2 = \SI{1.0}{\second}$ to illustrate the numerical 
penetration of the spherical obstacle by the \modif{plates, which is comparable to the finite element meshsize ($=1/64$). 
Note that the grey  area depicted in the top figure can bend as well since only its corners are fixed.
This effect can be reduced via  higher resolution and smaller penalty parameter $\varepsilon$}.}
\end{figure}

To simulate the essential effects
we consider a star-shaped configuration of a bilayer as shown in 
Figure~\ref{fig:microcap}. To start the encapsulation process we use
a Robin boundary condition on the entire boundary with external temperature
\[
\theta_{ext} = \SI{100.0}{\celsius}
\qquad
\text{on } \p_R\o=\p\o.
\]
The boundary is mechanically free but we fix the deformation 
at the vertices of one element in the mesh $\cT_h$ which contains
the midpoint of the center square to guarantee 
well-posedness of our method, i.e., uniquely defined iterates in
Algorithm~\ref{alg:alg_fully} (fully practical scheme). If the device center has the coordinate $(0,0)$, then the obstacle consists of $5$ spheres $B_i$, $i=0,...,4$, each of radius $\SI{0.24}{mm}$ and centered at $(0.28,0.28,0.25)$, $(0.72,0.28,0.25)$, $(0.28,0.72,0.25)$, $(0.72,0.72,0.25)$ and $(0.5,0.5,0.5)$ (all in $\SI{}{mm})$. At each vertex, the projection to the obstacle required by Algorithm~\ref{alg:alg_fully} is approximated by the projection to the closest sphere.
The obstacle models particles that are encapsulated by the deformed
bilayer plate.

We construct the mesh $\cT_h$ via~\modif{6} uniform quad refinements of
the initial coarse partition of the domain $\o$, indicated 
by the dashed lines in Figure~\ref{fig:microcap}. For our 
simulations we use the uniform time step \modif{$\tau = \SI{2.5e-4}{\second}$}
and the penalization parameter \modif{$\veps = \SI{5.0e-8}{mm^4  \megapascal^{-1}}$}.
\rhn{Since $T=1 s$ and $\mu_0=2\times10^3\megapascal$, this choice
satisfies \eqref{tau-epsilon}.}
Figure~\ref{f:microcap_sim} depicts the encapsulation process 
of~5 non-penetrable and rigid spherical particles with radii $\SI{0.24}{\mm}$. 
Penalization of the discrepancy between $y_h^{k+1}$ and
$s_h^k$ in Algorithm~\ref{alg:alg_fully} (fully practical scheme) does not prevent
penetration of the obstacle, an effect that depends on the size of 
the penalization parameter $\varepsilon$ \modif{ and the finite element meshsize and is discussed in Section~\ref{subsec:switch}}. We refer to the lower part
of Figure~\ref{f:microcap_sim} that illustrates this feature.
Reducing the size of $\varepsilon$ ameliorates this 
numerical artifact without making the algebraic system for $y_h^{k+1}$ of
Algorithm~\ref{alg:alg_fully} stiffer.

%%%%%%%%%%%%%%%%%%%%%%%%%%%%%%%%%%%%%%%%%%%%%%%%%%%%%%%%%%%%%%%%%%%%%%%%%%%%%%%%%%
\section{Conclusions}
%%%%%%%%%%%%%%%%%%%%%%%%%%%%%%%%%%%%%%%%%%%%%%%%%%%%%%%%%%%%%%%%%%%%%%%%%%%%%%%%%%

\rhn{This paper develops a reduced model for the thermal actuation of
bilayer plates consisting of a 4-th order nonlinear
PDE and a 2-nd order diffusion in the mid-surface of a thin plate.
The evolution is dictated by the latter whereas the former corresponds
to a quasi-stationary reaction of the plate.
We design a novel and effective finite element method for simulating
large 3D deformations of slender planar structures, including the
presence of obstacles. Several simulations illustrate the virtues of
this approach.

The planar structures must be macroscopically compliant to allow for
3D shape changes with relatively small external stimuli. We study thin
polymer materials composed of two layers which respond differently to
thermal actuation, resulting in out of plane curvature changes and
extremely large 3D deformations. 
The thin layer structures experience relatively small strain/stretch
while undergoing large rotation. The interface between the two layers,
or mid surface, is thus assumed to be perfectly bonded and inextensible,
which leads to isometric deformations.
The equations governing heat conduction and mechanical
deformation decouple, i.e., the thermal effect induces plate deformations, but
the latter do not affect heat conduction. This simplifies the solution process.

The algorithm takes advantage of several geometric properties valid
for isometries, and advances in time with a semi-implicit scheme; the
algebraic equations to be solved at each time step are linear. The
only restrictions among discretization parameters come from achieving
physically realistic dynamics in the presence of obstacles.
Simulations of significant practical applications examine the
feasibility of certain geometries and material parameters
of polymers in attaining dramatic shape changes. This methodology
could be used at the design stage, upon choosing material and
geometric parameters as well as external stimuli and simulating
shape reconfigurations prior to fabricating structures. This
methodology is general and can be extented for more complex material
models for active materials.}

%%%%%%%%%%%%%%%%%%%%%%%%%%%%%%%%%%%%%%%%%%%%%%%%%%%%%%%%%%%%%%%%%%%%%%%%%%%%%%%%%%

\bibliographystyle{acm}

%%%%%%%%%%%%%%%%%%%%%%%%%%%%%%%%%%%%%%%%%%%%%%%%%%%%%%%%%%%%%%%%%%%%%%%%%%%%%%%%%%
\end{document}